\documentclass[a4paper]{amsart}

\usepackage[T1]{fontenc}
\usepackage{textcomp}
\usepackage{lmodern}
\usepackage[latin1]{inputenc}
\usepackage[swedish,english]{babel}

\usepackage[pdftex]{color}
\usepackage{graphicx}
\usepackage{amsmath} 
\usepackage{amsthm}
\usepackage{amsfonts}
\usepackage{amssymb}
\usepackage{latexsym}

\newcommand{\eg}{\emph{e.g.}}
\newcommand{\ie}{\emph{i.e.}}
\newcommand{\ud}{\ \mathrm{d}}
\newcommand{\divop}{\mathrm{div}}
\newcommand{\rset}{\mathbb{R}}
\DeclareMathOperator*{\argmin}{argmin}

\begin{document}

% Document information
\title{Symplectic Reconstruction of Data for Heat and Wave Equations}

\thanks{Support by the Swedish Research Council grants 2002-6285 and
  2002-4961, and the European network HYKE, funded by the EC as
  contract HPRN-CT-2002-00282, is acknowledged.} 

\author{Jesper Carlsson}

\address{CSC, Numerical Analysis,
  Kungl. Tekniska H\"ogskolan,
  100 44 Stockholm, Sweden; 
}
\email{jesperc@kth.se}

%  65N21 - Partial differential equations, boundary value problems; Inverse problems
%  49L25 - Hamilton-Jacobi theories, including dynamic programming; Viscosity solutions


\begin{abstract}
  
  This report concerns the inverse problem of estimating a
  spacially dependent coefficient of a partial differential equation from
  observations of the solution at the boundary. 
  Such a problem can be formulated as an optimal control problem with
  the coefficient as the control variable and the solution
  as state variable. The heat or the
  wave equation is here considered as state equation. 
  It is well known that such inverse problems are ill-posed and need
  to be regularized. The powerful Hamilton-Jacobi theory is used to
  construct a simple and general method where the first step is to
  analytically regularize the 
  Hamiltonian; next its Hamiltonian system, a system of nonlinear
  partial differential equations, is solved with the Newton method
  and a sparse Jacobian. 

\end{abstract}

\maketitle
\tableofcontents

\section{Introduction}
%TODO: Short description about reconstruction, ill-posednes,
%regularizations and some references. 

In this paper we study the inverse problem to determine
a spacially dependent coefficient $\sigma$ of a partial differential equation
from partial knowledge of the forward solution $u$. In particular, we
seek the diffusion coefficient in the heat equation
and the wave speed coefficient in the wave equation.
Inverse problems arise in many applications such as inverse
scattering, impedance tomography and topology optimization, see \eg{} 
\cite{kunisch,bensoe,borcea,lions}, and share the
property that they are ill posed \ie{} given data $u$
there may not exist a corresponding coefficient $\sigma$, and if it
exists it may not be unique nor depend continuously on $u$. 
To be able to determine $\sigma$ the problem thus needs to be
regularized such that it becomes well posed.
The method used here to regularize and to solve the inverse problem 
is based on the work \cite{jesper,css,ss,mattias} where the inverse
problem is formulated as an optimal control problem and the
corresponding Hamilton-Jacobi equation is used to construct a
regularization, to obtain convergence results, and to finally solve the
regularized problem by using the method of characteristics \ie{}
to solve the corresponding Hamiltonian system.

The paper is stuctured as follows: In Section \ref{sec:opt} the
general theory of optimal control of partial differential equations
and Hamilton-Jacobi-Bellman is presented. In Section \ref{sec:heat}
the idea of how to optimally control the heat equation is discussed
together with numerical examples, and in Section
\ref{sec:wave} the control of the wave equation is treated.

\section{Optimal Control and Dynamic Programming}\label{sec:opt}
Consider a differential equation constrained minimization problem with
solution $\varphi:\Omega \times [0,T] \to V$, $\varphi =
\varphi(x,t)$ and control $\sigma:\Omega \times [0,T] \to B:= W$, $\sigma =
\sigma(x,t)$ for an open domain $\Omega$, some Hilbert space $V$ on $\Omega$,
and closed bounded set $B\subset \rset$:
\begin{equation}
  \label{eq:opt_problem}
  \begin{aligned}
    &\min_{\sigma: \Omega \times [0,T] \to B} 
    \int_0^T h(\varphi,\sigma) \ud t + g(\varphi^T),\\
    &\varphi_t = f(\varphi,\sigma),\\    
  \end{aligned}
\end{equation}
with $\varphi^T:=\varphi(\cdot,T)$ and given initial value
$\varphi^0=\varphi(\cdot,0)$. Here, $\varphi_t$ 
denotes the partial derivative with respect to time, $f: V \times
W \to V$ is the flux, and $h: V \times W \to \rset$, $g: V \to \rset$ are
given functions. 

This optimal control problem can be solved either directly using
constrained minimization or by dynamic programming.
The Lagrangian becomes 
\begin{equation*}
  {L}(\varphi,\lambda,\sigma):= \int_0^T
  \langle\lambda,f(\varphi,\sigma)-\varphi_t\rangle +
  h(\varphi,\sigma) \ud t,
\end{equation*}
with Lagrange multiplier $\lambda:\Omega \times [0,T] \to V$, and the
constrained minimization method is based on the Pontryagin method 
\begin{equation}
  \label{eq:pontryagin}
  \begin{aligned}
    \varphi_t & =f(\varphi,\sigma),\\
    \lambda_t & =-\langle \lambda,
    f_\varphi(\varphi,\sigma)\rangle + h_\varphi(\varphi,\sigma),\\
    \sigma(\cdot,t) & \in \argmin_{a:\Omega \to B} \{ \langle
    \lambda,f(\varphi,a)\rangle + h(\varphi,a) \}. 
  \end{aligned}
\end{equation}
with given initial value $\varphi^0$, final value
$\lambda^T:=\lambda(\cdot,T)=g_\varphi(\varphi^T)$, and where 
$f_\varphi$, $h_\varphi$ denotes the Gateaux derivatives with
respect to $\varphi$ and $\langle v,w \rangle$ is the duality pairing on
$V$, which reduces to the $L^2(\Omega)$ inner product if $v,w \in
L^2(\Omega)$.
For a differentiable Lagrangian that is convex in $\sigma$ the Pontryagin principle
coincides with the Lagrangian formulation for a constrained interior minimum
\begin{equation}
  \label{eq:KKT}
  \begin{aligned}
    \varphi_t & =f(\varphi,\sigma)\\
    \lambda_t & =-\langle \lambda,
    f_\varphi(\varphi,\sigma)\rangle + h_\varphi(\varphi,\sigma)\\ 
    0 & =\langle \lambda,
    f_\sigma(\varphi,\sigma)\rangle + h_\sigma(\varphi,\sigma),\\ 
    \sigma &\in B,\\
  \end{aligned}
\end{equation}
but in general \eqref{eq:pontryagin} and \eqref{eq:KKT} may have
different solutions $\varphi, \lambda, \sigma$ although both describe
necessary conditions for a minimizer to \eqref{eq:opt_problem}.
If an explicit minimizer in \eqref{eq:pontryagin} can be found
the Pontryagin principle gives additional information about
the control. Pontryagin's minimum principle can also be written as a
Hamiltonian system, see \cite{barron},
\begin{equation}
  \label{eq:HS}
  \begin{aligned}
    \varphi_t & ={H}_{\lambda} (\varphi, \lambda)\\
    \lambda_t & =-{H}_{\varphi} (\varphi, \lambda)\\
  \end{aligned}
\end{equation}
with $\varphi^0$ given,
$\lambda^T=g_\varphi(\varphi^T)$, and the Hamiltonian 
$H:V\times V\to \rset$ defined as
\begin{equation}
  \begin{aligned}
    \label{eq:hamiltonian}
    & H(\lambda,\varphi) &:=\min_{a: \Omega \to B}\{\langle
    \lambda,f(\varphi,a)\rangle + h(\varphi,a)\}.\\ 
  \end{aligned}
\end{equation}

The alternative dynamic programming method is based on the value
function $U:V\times [0,T] \to \rset$, 
\begin{equation*}
  \label{eq:value_function}
  \begin{aligned}
    &U(\phi,\tau):= \inf_{\sigma:\Omega\times[\tau,T]\to B} \bigg\{ \int_\tau^T
    h(\varphi,\sigma) \ud t +g(\varphi^T) \ \bigg |\
    \varphi_t = f(\varphi,\sigma), \ \varphi(\cdot,\tau)=\phi \in V
    \bigg\} \\ 
  \end{aligned}
\end{equation*}
which solves the nonlinear Hamilton-Jacobi-Bellman equation
\begin{equation}
  \label{eq:HJB}
  \begin{aligned}
    &\partial_t U(\phi,t) + {H} \big(U_\phi(\phi,t),\phi \big)=0,
    \quad U(\phi, T)=g(\phi),\\ 
  \end{aligned}
\end{equation}
with Hamiltonian defined as in \eqref{eq:hamiltonian}. Note that
solving the Hamiltonian system \eqref{eq:HS} is the method of
characteristics for the Hamilton-Jacobi equation \eqref{eq:HJB},
with $\lambda(x,t)=U_\varphi(\varphi(x,t),t)$.  
In general, the value function is however not
everywhere differentiable and the multiplier $\lambda$ becomes ill
defined in a classical sense. 

The Hamilton-Jacobi formulation \eqref{eq:HJB} has the advantages that
there is a complete well-posedness theory for Hamilton-Jacobi
equations, based on non-differential viscosity solutions, see \cite{cel}, and it finds a
global minimum. 
However, \eqref{eq:HJB} is not
computationally feasible for problems in high dimension, such as the
case where $\varphi$ is an approximation of a solution to a partial
differential equation. 
The Hamiltonian form  \eqref{eq:HS} has the advantage that it is
computationally feasible but the drawbacks are that it only focuses on
local minima and that the Hamiltonian
\eqref{eq:hamiltonian} in general only is Lipschitz continuous, even
if $f,g$ and $h$ are smooth, which means that the optimal control
depends discontinuously on $(\lambda,\varphi)$ and \eqref{eq:HS} becomes
undefined where the Hamiltonian is not differentiable.

In the following sections we will
%We will use a formulation of \eqref{eq:pontryagin} based on
%the Hamiltonian to regularize our problems and 
use a regularized version
of \eqref{eq:HS} to iteratively solve the nonlinear
constrained optimization problem \eqref{eq:opt_problem}.

% TODO: Need references about the generality (if exists) of the Pontryagin maximum
%principle for PDE constrained minimization. 

\section{Parameter Reconstruction for the Heat Equation }\label{sec:heat}
A distributed parameter reconstuction problem for the heat equation is
to find a heat 
conductivity (the control) \eg{} $\sigma : \bar \Omega \times [0,T] \to
[\sigma_-,\sigma_+]$, $\sigma=\sigma(x,t)$, $0<\sigma_-<\sigma_+$, and
a temperature distribution (the state) $u : \bar
\Omega \times [0,T] \to V$, $u=u(x,t)$ that satifies the heat equation 
\begin{equation} 
  \label{eq:constr_heat}
  \begin{aligned} 
    u_{t} &= \divop ( \sigma \nabla u ),
    & \text{ in } \Omega \times (0,T],\\
    \sigma \nabla u \cdot \bf{n} &= j, 
    & \text{ on } \partial\Omega \times (0,T],\\
    u &= 0,
    & \text{ on } \bar \Omega \times \{t=0\},
  \end{aligned}
\end{equation}
such that the error functional
\begin{equation} 
  \label{eq:obj_heat}
  \int_0^T \int_{\partial\Omega} (u-u^*)^2 \ud s \ud t,
\end{equation}
is minimized. The function $u^*=u^*(x,t)$ often represents
physical measurements contaminated by some noise, \eg{} $u^*(x,t) = u_{true}(x,t)
1 + w(x,t)$ where $w$ is a noise term and $u_{true}$ satisfies the above heat
equation for some unknown parameter $\sigma_{true}$, and in practice
the control is only spacially dependent, $\sigma_{true}=\sigma_{true}(x)$. 
The primary goal is thus to determine the unknown diffusion
coefficient $\sigma_{true}$ and the method to do so is to minimize the
objective functional 
\eqref{eq:obj_heat}. 

Inverse problems like \eqref{eq:constr_heat}, \eqref{eq:obj_heat} are in
general ill-posed due to one or more of the following reasons: 
\begin{enumerate}

\item There exists no minimizer $(u,\sigma)$, something that may occur
  with noisy data. Given unperturbed data
  $u^*$ corresponding to $\sigma_{true}$, it is evident that there
  exists a minimizer to \eqref{eq:constr_heat}, \eqref{eq:obj_heat}.

\item The minimizer is not unique, \eg{} although it may be possible to find an
  optimal state that minimizes \eqref{eq:obj_heat}, $u$ and $\sigma$
  may not be unique in $\Omega$.   

\item The solution $(u,\sigma)$, and particularly the control $\sigma$,
  depends discontinuously on data $u^*$. 

\end{enumerate}
A simple and common way to impose well-posedness to many inverse problems
is to add a Tikhonov regularization of the
form $\epsilon \| \sigma \|_{L^2(\Omega \times (0,T))}^2$ 
for  $\epsilon>0$, 
to the objective functional \eqref{eq:obj_heat}, see
\cite{engl,vogel,kunisch,lions}. Using the Pontryagin principle
presented in the previous section we will 
in Section \ref{sec:regularization_heat}
regularize the inverse problem \eqref{eq:constr_heat},
\eqref{eq:obj_heat} in a way that is comparable to a Tikhonov
regularization. 
% This regularization gives 
% an approximation of the unique value function and has stabilizing
% properties, 
% but in general the state and control may not be unique. 

% Some inverse problems can be exactly controllable, with a unique
% control, \ie{} for a 
% given $u^*$ there exist a unique $\sigma$ and $u=u^*$, see
% \cite{isakov,lions}, these problems are however usually still ill-posed in
% the sense of discontinuous dependency on data, but can with the addition
% of a Tikhonov regularization be completely well-posed.

Formulated as an optimal control problem the most
natural assumption on the control $\sigma$ is that it is dependent
on both time and space but as we will see in Section
\ref{sec:time-indep-control} it is also 
possible to let $\sigma = \sigma(x)$, $\sigma = \sigma(t)$, or even
let $\sigma$ be constant in time and space.  

% TODO: kunish paper

\subsection{The Hamiltonian System}\label{sec:hamilt_heat}
Following Section \ref{sec:opt} the Hamiltonian associated to the optimal 
control problem \eqref{eq:constr_heat} and \eqref{eq:obj_heat} is
\begin{equation}
  \label{eq:hamilt_heat}
  \begin{aligned}
  {H}(u,q,t) &:= \min_{\sigma : \Omega \to [\sigma_-,\sigma_+]} 
  \int_{\partial\Omega} (u-u^*)^2 \ud s + 
  \int_\Omega \divop ( \sigma \nabla u ) q \ud x\\
  &=\int_{\partial\Omega} (u-u^*)^2 + j q \ud s + 
  \min_{\sigma : \Omega \to [\sigma_-,\sigma_+]} 
   \int_\Omega - \sigma \nabla u \cdot \nabla q \ud x\\
   &=\int_{\partial\Omega} (u-u^*)^2 + j q \ud s - 
   \int_\Omega \underbrace{ 
     \max_{\sigma \in [\sigma_-,\sigma_+]} 
     \{ \sigma \nabla u \cdot \nabla q \} 
     }_{\mathfrak{h}(\nabla u \cdot \nabla q)} 
     \ud x.\\
 \end{aligned}
\end{equation}
and the Hamiltonian system, in strong form, then becomes
\begin{equation}
  \label{eq:hamilt_syst_heat}
  \begin{aligned}
    u_t &= \divop \big( \tilde \sigma \nabla
    u \big), & \text{ in } \Omega \times (0,T],\\
    \tilde \sigma \nabla u \cdot \bf{n} &= j, 
    & \text{ on } \partial\Omega \times (0,T],\\
    u &= 0, & \text{ on } \bar \Omega \times \{t=0\},\\
    -q_t &= \divop \big( \tilde \sigma \nabla
    q \big), & \text{ in } \Omega \times (0,T],\\
    \tilde \sigma \nabla q \cdot \bf{n} &=  2(u-u^*), 
    & \text{ on } \partial\Omega \times (0,T],\\
    q & = 0, & \text{ on }\Omega \times \{t=T\},\\
  \end{aligned}
\end{equation}
with 
\begin{equation}
  \label{eq:control}
  \tilde \sigma := \mathfrak{h}'(\nabla u \cdot \nabla q).
\end{equation}
It is here evident that the Hamiltonian only is Lipschitz continuous
and the control $\tilde \sigma$ is a bang-bang type control and
depends discontinuously on the
solutions $(u,q)$, see Figure \ref{fig:hamiltonian}. From the
optimality conditions \eqref{eq:KKT} an optimal solution has to satisfy 
$\nabla u \cdot \nabla q = 0$ and
\eqref{eq:hamilt_syst_heat} is thus undefined since $\mathfrak{h}'(0)$
is set valued, which calls
for a regularization.

\subsection{Regularization}\label{sec:regularization_heat}
\begin{figure}[htbp]
  \centering
  \includegraphics[width=1\textwidth]{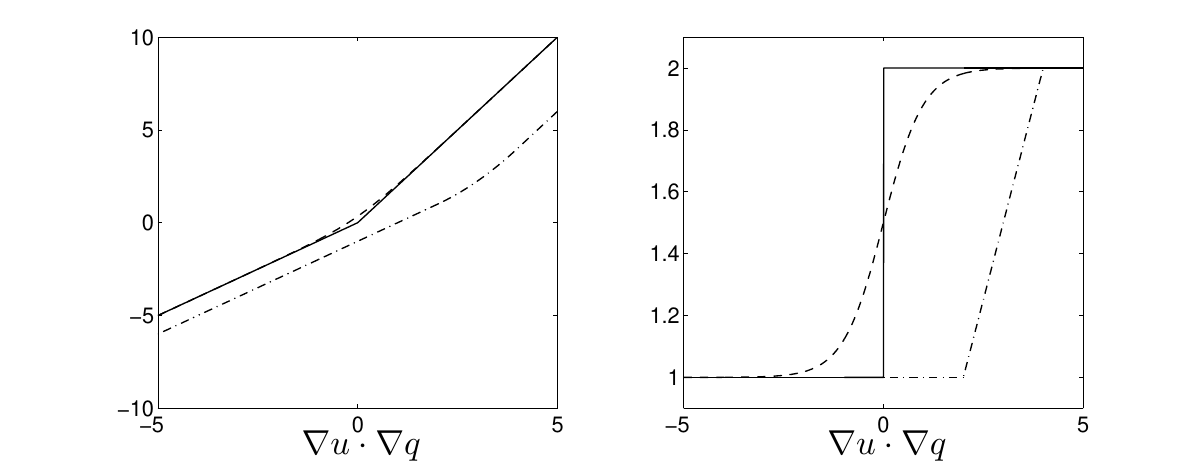}
  \caption{The functions $\mathfrak{h}$ (solid line),
    $\mathfrak{h}_\delta$ (dashed line),
    $\mathfrak{h}_{Tikhonov}$ (dash-dotted line) 
    to the left and their derivatives to the right.}
  \label{fig:hamiltonian}
\end{figure}

A simple regularization of the Hamiltonian system \eqref{eq:hamilt_syst_heat},
and consequently of the Hamiltonian \eqref{eq:hamilt_heat}, is to
approximate $\mathfrak{h}'$ with the
parabolic function
\begin{equation}
  \label{eq:reg}
  \mathfrak{h}'_\delta(\nabla u \cdot \nabla q) := 
  \underbrace{\frac{\sigma_+ + \sigma_-}{2}}_{\bar \sigma} + 
  \underbrace{\frac{\sigma_+ - \sigma_-}{2}}_{\hat \sigma} 
  \tanh(\frac{1}{\delta} \nabla u \cdot \nabla q),
\end{equation}
for some small $\delta>0$, see Figure \ref{fig:hamiltonian}. 
This regularization can be
compared with a classic Tikhonov regularization where
a small $L^2$-penalty of the control is added to the objective function
\eqref{eq:obj_heat}, \ie{} to minimize
\begin{equation} 
  \label{eq:obj_tikhonov}
  \int_0^T \int_{\partial\Omega} (u-u^*)^2 \ud s \ud t + 
  \delta \int_0^T \int_\Omega \sigma^2 \ud x \ud t.
\end{equation}
Minimizing \eqref{eq:obj_tikhonov} under the constraint
\eqref{eq:constr_heat} will lead to a $C^2$-Hamiltonian with 
\begin{equation*}
   H(u,q,t) = \int_{\partial\Omega} (u-u^*)^2 + j q \ud s - 
   \int_\Omega \underbrace{ 
     \max_{\sigma \in [\sigma_-,\sigma_+]} 
     \{ \sigma (\nabla u \cdot \nabla q -\delta \sigma) \} 
     }_{\mathfrak{h}_{Tikhonov}(\nabla u \cdot \nabla q)} 
     \ud x,
\end{equation*}
which can be seen in Figure \ref{fig:hamiltonian}. %Note that a $L^1$-penalty of the
%control will only result in a shift with respect to $\nabla u \cdot
%\nabla q$ in the Hamiltonian \eqref{eq:hamilt_heat} and does not act
%as a regularization but rather as a penalty.

\begin{figure}[htbp]
  \centering
  \includegraphics[width=0.6\textwidth]{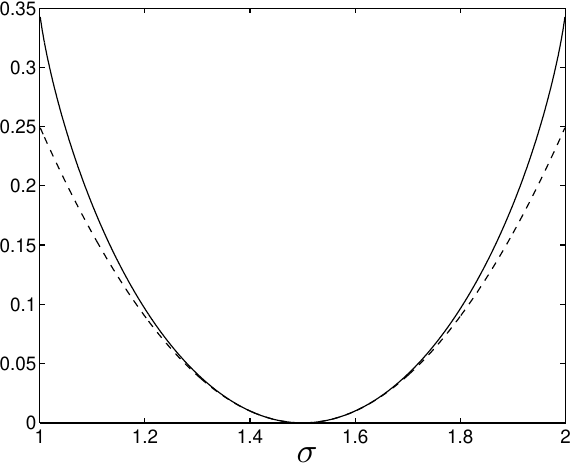}
  \caption{The function $\mathfrak{H}(\sigma)$ (solid line) compared to the
    $L^2$ penalty function $\delta (\sigma-\bar \sigma)^2$ (dashed line)
    for $\delta=1$, $\sigma_-=1$ and $\sigma_+=2$.}
  \label{fig:penalty}
\end{figure}

Another way to describe the simple regularization \eqref{eq:reg} 
is to see what kind of penalty on the objective function it
corresponds to.  We note that the
regularized Hamiltionian system can be written as 
\begin{equation*}
  \begin{aligned}
    \int_0^T \bigg( \int_\Omega -u_{t} v - 
    \mathfrak{h}_\delta'(\nabla u \cdot \nabla q) \nabla u \cdot \nabla v \ud x + 
    \int_{\partial\Omega} j v \ud s \bigg) \ud t &= 0, & \forall v
    \in V,\\
    \int_0^T \bigg( \int_\Omega q_{t} v - 
    \mathfrak{h}_\delta'(\nabla u \cdot \nabla q) \nabla q \cdot \nabla v \ud x + 
    \int_{\partial\Omega} 2(u-u^*) v \ud s \bigg) \ud t &= 0, &
    \forall v \in V,
  \end{aligned}
\end{equation*}
or by a redefinition of $\sigma$
% such that $\sigma:[0,T]\times \Omega
%\to W$ for some Hilbert space $W$
\begin{equation}
  \label{eq:hamilt_syst_alt}
  \begin{aligned}
    \int_0^T \bigg( \int_\Omega -u_{t} v - 
    \sigma \nabla u \cdot \nabla v \ud x + 
    \int_{\partial\Omega} j v \ud s \bigg) \ud t &= 0, & \forall v
    \in V,\\
    \int_0^T \bigg( \int_\Omega q_{t} v - 
    \sigma \nabla q \cdot \nabla v \ud x + 
    \int_{\partial\Omega} 2(u-u^*) v \ud s \bigg) \ud t &= 0, &
    \forall v \in V,\\
    \int_0^T \int_\Omega \Big( 
    \sigma - \mathfrak{h}_\delta'(\nabla u \cdot \nabla q) 
    \Big) v \ud x \ud t&= 0, & \forall v \in W,
  \end{aligned}
\end{equation}
where $\sigma:[0,T]\times \Omega \to W$ for some Hilbert space $W$.
Let $\mathfrak{H}$ be the primitive function of the inverse function
of $\mathfrak{h}_\delta'$ \ie{}
\begin{equation*}
  % \mathfrak{H}(\sigma) := \delta
  % \Big(
  % \sigma - \bar \sigma
  % \Big)
  % \text{tanh}^{-1}
  % \Big(
  % \frac{\sigma - \bar \sigma}{\hat \sigma}
  % \Big) + 
  % \frac{\delta \hat \sigma}{2} \ln
  % \bigg( 1 -
  % \Big(
  % \frac{\sigma - \bar \sigma}{\hat \sigma}
  % \Big)^2
  % \bigg),
  % \mathfrak{H}(\sigma) :=
  % \frac{\delta}{2}\Big(
  % (1+\frac{\sigma-\bar \sigma}{\hat \sigma}) 
  % \ln(1+\frac{\sigma-\bar \sigma}{\hat \sigma}) +
  % (1-\frac{\sigma-\bar \sigma}{\hat \sigma}) 
  % \ln(1-\frac{\sigma-\bar \sigma}{\hat \sigma})
  % \Big),
  \mathfrak{H}(\sigma) :=
  \frac{\delta}{2\hat \sigma}\bigg(
  (\sigma-\sigma_-)
  \ln\Big(\frac{\sigma-\sigma_-}{\hat \sigma}\Big) +
  (\sigma_+-\sigma) 
  \ln\Big(\frac{\sigma_+-\sigma}{\hat \sigma}\Big)
  \bigg),
\end{equation*}
then it is evident that \eqref{eq:hamilt_syst_alt} can be seen as the first
order optimality conditions for the problem to minimize
\begin{equation*} 
  \int_0^T \int_{\partial\Omega} (u-u^*)^2 \ud s \ud t + 
  \int_0^T \int_\Omega  \mathfrak{H}(\sigma) \ud x \ud t,
\end{equation*}
under the constraint \eqref{eq:constr_heat}. In Figure
\ref{fig:penalty} the function $\mathfrak{H}(\sigma)$ is compared with
a Tikhonov regularization of the form $\delta (\sigma-\bar \sigma)^2$.

It is often beneficial to prevent spacial oscillations of the
coefficient by adding a penalty on the $L^2$-norm of the
gradient of the coefficient, \ie{} 
$\epsilon \| \nabla \sigma \|_{L^2(\Omega \times (0,T))}^2$, for $\epsilon>0$,
to the objective function \eqref{eq:obj_heat}. 
For such a penalty the minimization in the corresponding
Hamiltonian
\begin{equation}
  \label{eq:hamiltonian_spacial_penalty}
  H(u,q,t):=\min_{\sigma : \Omega \to [\sigma_-,\sigma_+]} 
  \int_{\partial\Omega} (u-u^*)^2 \ud s + 
  \int_\Omega \divop ( \sigma \nabla u ) q + \epsilon |\nabla \sigma|^2 \ud x,
\end{equation}
can not be done explicitly, and instead taking the first variation in
$\sigma$ would give the system
\begin{equation*}
  \begin{aligned}
    u_t &= \divop \big( \sigma \nabla u \big),\\ 
    -q_t &= \divop \big( \sigma \nabla q \big),\\ 
    2\epsilon \Delta \sigma &= - \nabla u \cdot \nabla q, \\
    \sigma &\in [\sigma_-,\sigma_+].
  \end{aligned}
\end{equation*}
which corresponds to the usual first order optimality conditions for
the Lagrangian. How to treat different penalties on the control in an
optimal control setting is discussed in Section \ref{sec:penalty}.

% If we let
% \begin{equation*}
%   \mathfrak{H}(\sigma) :=   
%   - \frac{\delta \hat \sigma}{2} \ln
%   \bigg( 1 -
%   \Big(
%   \frac{\sigma - \bar \sigma}{\hat \sigma}
%   \Big)^2
%   \bigg),
% \end{equation*}
% then
% \begin{equation*}
%    \mathfrak{H}'(\sigma) = \mathfrak{h}'^{-1}(\sigma) = 
%    \delta \frac{\sigma-\bar\sigma}{\hat\sigma}
%    \bigg( 
%    1- \Big( \frac{\sigma-\bar\sigma}{\hat\sigma} \Big)^2 
%    \bigg)^{-1}
% \end{equation*}
% and
% \begin{equation*}
%   \mathfrak{h}'(\nabla u \cdot \nabla q) = \bar\sigma - 
%   \frac{\delta\hat\sigma}{2 \nabla u \cdot \nabla q} +
%   \sqrt{
%     \Big(
%     \frac{\delta\hat\sigma}{2 \nabla u \cdot \nabla q}
%     \Big)^2 + \nabla u \cdot \nabla q
%   }
% \end{equation*}

\subsection{Time Independent Control} \label{sec:time-indep-control} 
To study the case when the control $\sigma$ is
independent of time we first assume that it not only is independent of
time but also depends on an auxilliary variable $z$,
\ie{} $\sigma : \bar \Omega \times [0,\tilde T] \to
[\sigma_-,\sigma_+]$, $\sigma=\sigma(x,z)$. For a moment we also
assume that $u : \bar \Omega \times [0,T] \times [0,\tilde T] \to V$,
$u=u(x,t,z)$, but with the same measurements as in
\eqref{eq:obj_heat}.  If we treat $z$ as the time and $t$ as a spacial  
variable we can define the optimal control problem
\begin{equation}
  \label{eq:obj_heat_z}
  \min_{\sigma : \bar \Omega \times [0,\tilde T] \to [\sigma_-,\sigma_+]} 
  \frac{1}{\tilde T} \int_0^{\tilde T} \int_0^T \int_{\partial\Omega}
  (u-u^*)^2 \ud s \ud t \ud z,
\end{equation}
where the state $u$ satisfies the partial differential equation
\begin{equation} 
  \label{eq:constr_heat_z}
  \begin{aligned} 
    u_{z} &=\frac{1}{\tilde T} \Big( \divop ( \sigma \nabla u ) -u_t \Big),
    & \text{ in } \Omega \times (0,T) \times (0,\tilde T],\\
    \sigma \nabla u \cdot \bf{n} &= j, 
    & \text{ on } \partial\Omega \times (0,T) \times (0,\tilde T],\\
    u &= 0,
    & \text{ on } \bar \Omega \times \{t=0\} \times (0,\tilde T],\\
    u &= u_{0},
    & \text{ on } \bar \Omega \times (0,T) \times \{z=0\}.
  \end{aligned}
\end{equation}
for some arbitrary initial condition $u(x,t,0)=u_{0}$. 
%TODO It is easy to see that for a given $\sigma$, constant in $z$, the solution
%$u(x,t,\tilde T)$ to \eqref{eq:constr_heat_z} would simply approach
%the solution to \eqref{eq:constr_heat} as $\tilde T \to \infty$ since
%the boundary conditions does not depend on $z$. 

The Hamiltonian for \eqref{eq:obj_heat_z}, \eqref{eq:constr_heat_z} is
\begin{equation}
  \label{eq:time_indep_hamilt}
  \begin{aligned}
    {H}(u,q,z) :=& \min_{\sigma : \Omega \to [\sigma_-,\sigma_+]} 
    \frac{1}{\tilde T} \int_0^T \int_{\partial\Omega} (u-u^*)^2 \ud s \ud t\\
    &+\frac{1}{\tilde T} \int_0^T \int_\Omega 
    \Big( \divop ( \sigma \nabla u ) - u_t \Big)
    q \ud x \ud t\\
    =& \frac{1}{\tilde T} \int_0^T \int_{\partial\Omega} (u-u^*)^2 + j q \ud s \ud t - 
    \frac{1}{\tilde T} \int_0^T \int_\Omega u_t q \ud x \ud t\\
    & -\frac{1}{\tilde T} \int_\Omega \underbrace{
      \max_{\sigma \in [\sigma_-,\sigma_+]} \bigg\{
      \sigma\int_0^T \nabla u \cdot \nabla q  \ud t \bigg\}
    }_
    {\mathfrak{h} \big(\int_0^T \nabla u \cdot \nabla q
      \ud t \big)} \ud x,\\
  \end{aligned}
\end{equation}
and the Hamiltonian system is given by
\begin{equation}
  \label{eq:hamilt_syst_heat_z}
  \begin{aligned}
    u_z &= \frac{1}{\tilde T} \Big( \divop ( \mathfrak{h}' \nabla u ) - u_t \Big), 
    & \text{ in } \Omega \times (0,T) \times (0,\tilde T],\\
     \mathfrak{h}' \nabla u \cdot \bf{n} &= j, 
    & \text{ on } \partial\Omega \times (0,T) \times (0,\tilde T],\\
    u &= 0,
    & \text{ on } \bar \Omega \times \{t=0\} \times (0,\tilde T],\\
    u &= u_{0},
    & \text{ on } \bar \Omega \times (0,T) \times \{z=0\},\\
    -q_z &= \frac{1}{\tilde T} \Big( \divop ( \mathfrak{h}' \nabla q ) + q_t \Big), 
    & \text{ in } \Omega \times (0,T) \times (0,\tilde T],\\
     \mathfrak{h}' \nabla q \cdot \bf{n} &= 2(u-u^*), 
    & \text{ on } \partial\Omega \times (0,T) \times (0,\tilde T],\\
    q &= 0,
    & \text{ on } \bar \Omega \times \{t=T\} \times (0,\tilde T],\\
    q &= 0,
    & \text{ on } \bar \Omega \times (0,T) \times \{z=\tilde T\}.\\
  \end{aligned}
\end{equation}
Under the assumption that the solutions $u$ and $q$ in
\eqref{eq:hamilt_syst_heat_z} are asymptotically stationary as
$\tilde T \to \infty$, the Hamiltonian system for the problem
\eqref{eq:constr_heat}, \eqref{eq:obj_heat}, with a time-independent
control, is given by \eqref{eq:hamilt_syst_heat} and
\begin{equation}
  \label{eq:control_time_indep}
  \tilde \sigma := \mathfrak{h}'
  \bigg( \int_0^T \nabla u \cdot \nabla q \ud t \bigg).
\end{equation}
Similarly, the case of a space independent coefficient $\sigma=\sigma(t)$
will lead to 
\begin{equation*}
  \tilde \sigma := \mathfrak{h}'
  \bigg( \frac{1}{|\Omega|}\int_\Omega \nabla u \cdot \nabla q \ud x \bigg),
\end{equation*}
and for the case where $\sigma$ is constant
\begin{equation*}
  \tilde \sigma := \mathfrak{h}'
  \bigg( \frac{1}{|\Omega|} \int_0^T \int_\Omega \nabla u \cdot \nabla
  q \ud x \ud t \bigg).
\end{equation*}

\subsection{Penalty on the Control}\label{sec:penalty}
If we want to reconstruct a time independent control it can be
beneficial to put a penalty on $\sigma_t$, \ie{} we want to minimize the
objective functional
\begin{equation} 
  F(u,\sigma_t) := \int_0^T \int_{\partial\Omega} (u-u^*)^2 \ud s \ud t +  
  \varepsilon \int_0^T \int_{\Omega} \sigma_t^2 \ud x \ud t,
\end{equation}
under the usual constraint \eqref{eq:constr_heat}. To do this the
optimal control problem has to be reformulated such that $\sigma$ is
a state variable and the control is defined as $z:=\sigma_t(x,t)$, $z
: \bar \Omega \times [0,T] \to [z_-,z_+]$.
The optimal control problem is thus to find a control $z$ and state
variables $u$ and $\sigma$ such that $F(u,z)$ is minimized and the
system 
\begin{equation*} 
  \begin{aligned} 
    u_{t} &= \divop ( \sigma \nabla u ),
    & \text{ in } \Omega \times (0,T],\\
    \sigma_{t} &= z 
    &\text{ in } \Omega \times (0,T],\\
    \sigma \nabla u \cdot \bf{n} &= j, 
    & \text{ on } \partial\Omega \times (0,T],\\
    u &= 0,
    & \text{ on } \bar \Omega \times \{t=0\},\\
    \sigma &= \sigma_0 > 0,
    & \text{ on } \bar \Omega \times \{t=0\}.\\
  \end{aligned}
\end{equation*}
is satisfied.
The Hamiltonian becomes
\begin{equation*}
  \begin{aligned}
  {H}(u,q,\sigma,\lambda,t) :=& \min_{z : \Omega \to [z_-,z_+]} 
  \int_{\partial\Omega} (u-u^*)^2 \ud s +
  \int_\Omega \divop ( \sigma \nabla u ) q + z\lambda + \varepsilon z^2\ud x\\
  =&\int_{\partial\Omega} (u-u^*)^2 + j q \ud s - 
  \int_\Omega \sigma \nabla u \cdot \nabla q \ud x\\
  &+\int_\Omega \underbrace{ 
    \min_{z : \Omega \to [z_-,z_+]} \{ z (\varepsilon z + \lambda) \} 
  }_{\mathfrak{h}(\lambda)} \ud x,\\
   \end{aligned}
\end{equation*}
and the corresponding Hamiltonian system is
\begin{equation*}
  \begin{aligned}
    u_t &= \divop \big( \sigma \nabla
    u \big), & \text{ in } \Omega \times (0,T],\\
    \sigma_{t} &= \mathfrak{h}'(\lambda)
    &\text{ in } \Omega \times (0,T],\\
    \sigma \nabla u \cdot \bf{n} &= j, 
    & \text{ on } \partial\Omega \times (0,T],\\
    u &= 0, & \text{ on } \bar \Omega \times \{t=0\},\\ 
    \sigma &= \sigma_0 > 0,
    & \text{ on } \bar \Omega \times \{t=0\}.\\
    -q_t &= \divop \big( \sigma \nabla
    q \big), & \text{ in } \Omega \times (0,T],\\
    -\lambda_{t} &= - \nabla u \cdot \nabla q
    &\text{ in } \Omega \times (0,T],\\
    \sigma \nabla q \cdot \bf{n} &=  2(u-u^*), 
    & \text{ on } \partial\Omega \times (0,T],\\
    q & = 0, & \text{ on }\Omega \times \{t=T\},\\
    \lambda & = 0, & \text{ on }\Omega \times \{t=T\},\\
  \end{aligned}
\end{equation*}
which is equivalent to \eqref{eq:hamilt_syst_heat} with 
\begin{equation*}
  \tilde \sigma := \sigma_0 + 
  \int_0^t \mathfrak{h}' \bigg( 
  \int_y^T -(\nabla u \cdot \nabla q)(x,z) \ud z 
  \bigg) \ud y.
\end{equation*}
Note, since we no longer have a constraint $\sigma>0$, the bound $z_-$
has to be carefully chosen to ensure well-posedness of the forward problem.

In a similar fashion as for penalizing temporal variations of the
control it is also possible to penalize spacial variations, as was
briefly mentioned in Section \ref{sec:regularization_heat}, where the
objective was to minimize $F(u,|\nabla \sigma|)$ under the constraint
\eqref{eq:constr_heat}, which leads to the Hamiltonian
\eqref{eq:hamiltonian_spacial_penalty}. To be able to explicitly find
the minimum in the Hamiltonian we once again let $\sigma$
act as a state variable, introduce the control $z$ and the dynamics 
\begin{equation}
  \label{eq:dynamics}
  \begin{aligned}
    \sigma_t & = \frac{z-|\nabla\sigma|^2}{\gamma}, & \text{ in } \Omega
    \times (0,T],\\
    \sigma & = \sigma_0 > 0, & \text{ in } \Omega
    \times \{t=0\},\\ 
  \end{aligned}
\end{equation}
for $\gamma>0$.
The slightly perturbed control problem is now to minimize the objective
function $F(u,z)$ such that \eqref{eq:constr_heat} and
\eqref{eq:dynamics} holds, which leads to the Hamiltonian
\begin{equation*}
  \begin{aligned}
    {H}(u,q,\sigma,\lambda,t) :=& \min_{z : \Omega \to [z_-,z_+]} 
    \int_{\partial\Omega} (u-u^*)^2 \ud s +
    \int_\Omega \divop ( \sigma \nabla u ) q + \lambda\frac{z-|\nabla\sigma|^2}{\gamma} + \varepsilon z\ud x\\
    =&\int_{\partial\Omega} (u-u^*)^2 + j q \ud s - 
    \int_\Omega \sigma \nabla u \cdot \nabla q + \lambda\frac{|\nabla\sigma|^2}{\gamma} \ud x\\
    &+\int_\Omega \underbrace{ 
      \min_{z : \Omega \to [z_-,z_+]} \{ z (\varepsilon + \frac{\lambda}{\gamma}) \} 
    }_{\mathfrak{h}(\lambda)} \ud x,\\
  \end{aligned}
\end{equation*}
and the Hamiltonian system
\begin{equation*}
  \begin{aligned}
    u_t &= \divop \big( \sigma \nabla u \big),\\ 
    \sigma_t &= \mathfrak{h}'(\lambda)-\frac{|\nabla\sigma|^2}{\gamma},\\
    -q_t &= \divop \big( \sigma \nabla q \big),\\ 
    -\lambda_t &= \nabla u \cdot \nabla q - 2\lambda \frac{\Delta \sigma}{\gamma}.\\ 
  \end{aligned}
\end{equation*}

\subsection{Numerical Approximation and Symplectic Methods} \label{sec:symplectic_heat}
Let $\bar V \subset V:=H^1(\Omega)$ be the finite element subspace of piecewise
linear functions defined on a triangulation of $\Omega$, which implies
that our optimal control problems in the previous sections are
approximated by optimal control problems for ordinary differential
equations. 
We also let the functions $\mathfrak{h}_\delta$,$H^\delta$ and
$h_\delta$ denote the regularized counterparts to $\mathfrak{h}$,$H$ and
$h$. The regularized version of
$\mathfrak{h}$ is given by \eqref{eq:reg} from which the definition of $H^\delta$
follows. The regularized function $h_\delta$ can be derived from $H^\delta$
by $h_\delta:=H^\delta-\langle\lambda, H_\lambda^\delta \rangle$ and a
regularized version of $f$ can be defined as $f_\delta:=H_\lambda^\delta$. 

Now, introduce the uniform partition $\{t_i=ki\}_{i=0}^N$, $k=T/N$ of the
time interval $[0,T]$, and the 
corresponding finite element approximations at each time step
$\varphi_n:=\varphi(t_n), \lambda_n:=\lambda(t_n)$. 
Also define a discrete regularized
version $\bar U: \bar V \times [0,T] \to \rset$ of the value function
\eqref{eq:value_function}, 
\begin{equation*}
  \label{eq:value_disc}
  \bar U(\phi,t_m):= \min_{\varphi_m=\phi} \bigg\{
  g(\varphi_N) + k \sum_{n=m}^{N-1}
  h_\delta(\varphi_n,\lambda_{n+1})\bigg\},
\end{equation*}
where $\varphi_n$ and $\lambda_n$ satisfy a symplectic scheme, \eg{}
the symplectic forward Euler method
\begin{equation}
  \label{eq:sympl_forw_euler}
  \begin{aligned}
    \varphi_{n+1}-\varphi_n &= k {H}_\lambda^\delta ( \varphi_{n}, \lambda_{n+1} ), 
    && \text{ for } n=m,\ldots,N-1 \text{ given }\varphi_m = \phi,\\
    \lambda_{n}-\lambda_{n+1} &= k {H}_\varphi^\delta ( \varphi_{n}, \lambda_{n+1} ), 
    &&  \text{ for } n=m,\ldots,N-1 \text{ given } \lambda_N = g_\varphi(\varphi_N).\\
  \end{aligned}
\end{equation}
Symplecticity here means that $\bar U_\varphi(\varphi_n,t_n) = \lambda_n$, 
\ie{} the gradient of the discrete value function coincides with the
discrete dual $\lambda_n$, and given that
$|H-H^\delta|=\mathcal{O}(\delta)$ it can be shown that for symplectic
one-step schemes 
\begin{equation*}
  \bigg| U(\varphi_0,t_0) - 
  g(\varphi_N) - k \sum_{n=m}^{N-1}
  h_\delta(\varphi_n,\lambda_{n+1}) \bigg| = \mathcal{O}(k),
\end{equation*}
for $\delta \sim k$, see \cite{ss}.
It is thus essential to use a symplectic time discretization of the
regularized Hamiltonian system
\begin{equation*}
  \begin{aligned}
    \varphi_t &= {H}_\lambda^\delta ( \varphi, \lambda), \\
    \lambda_t &= - {H}_\varphi^\delta ( \varphi, \lambda ), 
  \end{aligned}
\end{equation*}
in order to have convergence in the value function.
%It can be shown that under certain assumptions this approximated value
%function solves a modified Hamilton-Jacobi equation and approximates
%the true value function $U: \bar V \times [0,T] \to \rset$ (on a finite element
%space) with an error 
%of order $k$ independently of the spacial discretization, see
%\cite{mattias} and \cite{css}. 

%To appropriately discretice the Hamiltonian system \eqref{eq:HS} it is
%essential to use a symplectic time discretization for $\varphi_n,
%\lambda_n$ which in this case means
%that 
Some examples of other symplectic schemes are the the backward Euler method 
\begin{equation}
  \label{eq:sympl_backw_euler}
  \begin{aligned}
    \varphi_{n+1}-\varphi_n &= k {H}_\lambda^\delta ( \varphi_{n+1}, \lambda_{n} ), 
    && \text{ for } n=0,\ldots,N-1 \text{ given }\varphi_0,\\
    \lambda_{n}-\lambda_{n+1} &= k {H}_\varphi^\delta ( \varphi_{n+1}, \lambda_{n} ), 
    &&  \text{ for } n=0,\ldots,N-1 \text{ given } \lambda_N,\\
  \end{aligned}
\end{equation}
and the implicit midpoint method
\begin{equation}
  \label{eq:sympl_midpoint}
  \begin{aligned}
    \varphi_{n+1}-\varphi_n &= k {H}_\lambda^\delta 
    \bigg( \frac{\varphi_n+\varphi_{n+1}}{2}, \frac{\lambda_n+\lambda_{n+1}}{2} \bigg),
    && \text{ for } n=0,\ldots,N-1 \text{ given }\varphi_0,\\
    \lambda_{n}-\lambda_{n+1} &= k {H}_\varphi^\delta 
    \bigg( \frac{\varphi_n+\varphi_{n+1}}{2}, \frac{\lambda_n+\lambda_{n+1}}{2} \bigg),
    && \text{ for } n=0,\ldots,N-1 \text{ given }\lambda_N.\\
  \end{aligned}
\end{equation}
See \cite{wanner} for a thorough description of symplectic methods.

% Note that \eqref{eq:sympl_backw_euler} is the stationary point to the
% discrete Lagrangian
% \begin{equation*}
%   L:=\sum_{n=0}^{N-1} \langle \varphi_{n}-\varphi_{n+1}, \lambda_n \rangle + kH(\varphi_{n+1},\lambda_n),
% \end{equation*}

\subsection{The Newton Method}\label{sec:newton_heat}
To solve the coupled nonlinear symplectic schemes
\eqref{eq:sympl_forw_euler}-\eqref{eq:sympl_midpoint} above, it is tempting
to propose fix-point schemes that partly removes the coupling between the forward and
bacward equation, \eg{} by iterating separately in $\varphi$ and
$\lambda$. Such methods has the advantage that existing
partial differential equation
solvers can be used to efficiently solve the forward and backward
problems in each iteration, but the disadvantage is that the
convergence to an optimal solution tends to be slow, and also
dependent on the discretisation. A more suitable strategy is to use
information of the Hessian of $H^\delta$; \eg{} Quasi-Newton methods, or since
the Hessian in our case can be found explicitly and is sparse, the
Newton method itself.  

For the Hamiltonian system \eqref{eq:hamilt_syst_heat} with 
$\tilde \sigma:=\mathfrak{h}_\delta'$ given by \eqref{eq:reg} the symplectic
backward Euler can be written as
% \begin{equation}
%   \label{eq:hamilt_syst_heat_disc}
%   \begin{aligned}
%     \int_\Omega (u_{n+1}-u_n)v \ud x =& -k \int_\Omega
%     \mathfrak{h}'(\nabla u_{n+1} \cdot \nabla q_{n}) \nabla
%     u_{n+1} \cdot \nabla v \ud x + \\
%     &k\int_{\partial\Omega} j_{n+1}v \ud s,\\
%     \int_\Omega (q_{n}-q_{n+1})w \ud x =& -k \int_\Omega
%     \mathfrak{h}'(\nabla u_{n+1} \cdot \nabla q_{n}) \nabla
%     q_{n+1} \cdot \nabla w \ud x +\\
%     &k\int_{\partial\Omega} 2(u_{n+1}-u_{n+1}^*)w \ud s,\\
%   \end{aligned}
% \end{equation}
\begin{equation*}
    F_n(w) = 0, \quad G_n(w) = 0, \quad n=0,\ldots,N-1, \quad \forall
    w \in \bar V
\end{equation*}
where
\begin{equation}
  \label{eq:fg}
  \begin{aligned}
    F_n (w):= &\int_\Omega (u_{n+1}-u_n) w + k \mathfrak{h}_\delta'(\nabla
    u_{n+1} \cdot \nabla q_{n}) \nabla u_{n+1} \cdot \nabla w \ud x\\
    &-\int_{\partial\Omega} k j_{n+1} w \ud s,\\   
    G_n (w):= &\int_\Omega (q_{n}-q_{n+1}) w + k \mathfrak{h}_\delta'(\nabla
    u_{n+1} \cdot \nabla q_{n}) \nabla q_n \cdot \nabla w \ud x\\
    &-\int_{\partial\Omega} 2k(u_{n+1}-u^*_{n+1}) w \ud s,\\
  \end{aligned}
\end{equation}
and $u_0=q_N=0$. 
Given an initial guess $u[0]$, $q[0]$ the (damped) Newton method yields that 
\begin{equation*}
  \begin{aligned}
    u[i+1] = u[i] - \alpha \hat u,\\
    q[i+1] = q[i] - \alpha \hat q,\\
  \end{aligned}
\end{equation*}
where $\alpha\in(0,1]$ and, for each iteration, the updates $\hat
u$ and $\hat q$ solve a linear system of the form
\begin{equation}
  \label{eq:saddle_syst}
    \left(
    \begin{array}{cc}
      K_{11} & K_{12}\\
      K_{21} & K_{11}^T
    \end{array}
  \right)
  \left(
    \begin{array}{c}
      \hat u\\
      \hat q
    \end{array}
  \right)
  = 
  \left(
    \begin{array}{c}
      f\\
      g
    \end{array}
  \right),
\end{equation}
where
\begin{equation*}
  \begin{aligned}
    \hat u = &
    \left(
      \begin{array}{ccccc}
        \hat u_1 & \ldots & \hat u_N
      \end{array}
    \right)^T, 
    & \hat q = &
    \left(
      \begin{array}{ccccc}
        \hat q_0 & \ldots & \hat q_{N-1}
      \end{array}
    \right)^T,\\
    f = &
    \left(
      \begin{array}{ccccc}
        F_0 & \ldots & F_{N-1}
      \end{array}
    \right)^T,
    & g = &
    \left(
      \begin{array}{ccccc}
        G_0 & \ldots & G_{N-1}
      \end{array}
    \right)^T.\\
\end{aligned}
\end{equation*}
The matrix $K_{11}$ is a bi-diagonal block matrix with $M+S_i$ for $i=0,\ldots,N-1$
on the diagonal and $-M$ on the sub-diagonal, where $M$ denotes the mass matrix
\begin{equation*}
  \int_\Omega w\bar w \ud x,
\end{equation*}
and
\begin{equation*}
  \begin{aligned}
    S_n :=&  \int_\Omega k \mathfrak{h}_\delta''(\nabla
    u_{n+1} \cdot \nabla q_{n}) \nabla q_n \cdot \nabla w \
    \nabla u_{n+1} \cdot \nabla \bar w \ud x\\
    &+\int_\Omega k \mathfrak{h}_\delta'(\nabla
    u_{n+1} \cdot \nabla q_{n}) \nabla w \cdot \nabla \bar w \ud x.\\
  \end{aligned}
\end{equation*}
for $w, \bar w \in \bar V$.
The matrices $K_{12}$, $K_{21}$ are symmetric block-diagonal matrices with
\begin{equation*}
  \int_\Omega k \mathfrak{h}_\delta''(\nabla
    u_{n+1} \cdot \nabla q_{n}) \nabla u_{n+1} \cdot \nabla w \
    \nabla u_{n+1} \cdot \nabla \bar w \ud x,
\end{equation*}
and
\begin{equation*}
  \int_\Omega k \mathfrak{h}_\delta''(\nabla
  u_{n+1} \cdot \nabla q_{n}) \nabla q_n \cdot \nabla w \
  \nabla q_n \cdot \nabla \bar w \ud x -
  \int_{\partial\Omega} 2 k \bar w w \ud s,
\end{equation*}
for $n=0,\ldots,N-1$ on the the diagonals, respectively.

If we repartition the block $2 \times 2$ linear system
\eqref{eq:saddle_syst} to
\begin{equation}
  \label{eq:saddle_syst2}
    \left(
    \begin{array}{cc}
      K_{21} & K_{11}^T\\
      K_{11} & K_{12}
    \end{array}
  \right)
  \left(
    \begin{array}{c}
      \hat u\\
      \hat q
    \end{array}
  \right)
  = 
  \left(
    \begin{array}{c}
      g\\
      f
    \end{array}
  \right),
\end{equation}
we see that it is a generalized saddle point system
\cite{saddle-point} with symmetric matrices $K_{21},
K_{12}$, and $K_{11}^T \neq 0$, $K_{21}\neq 0$. However, unlike saddle
point problems arising from \eg{} the steady-state Navier-Stokes
equations or from the Karush-Kuhn-Tucker optimality conditions for
equality constrained minimization problems, both $K_{12}$ and $K_{21}$ may
here be indefinite and singular.

Since \eqref{eq:saddle_syst} and \eqref{eq:saddle_syst2} are
increasingly ill-conditioned with 
respect to reduction in mesh size, step size and
regularization, the success of iterative algorithms like Krylov
sub-space methods will depend heavily on the choice of preconditioner.
Standard algebraic
preconditioners like incomplete LU-factorization are often unsuitable
for saddle-point problems due to the indefiniteness and lack of
diagonal dominance, so the preconditioner must be tailored
for the specific problem at hand.
One popular approach for PDE-constrained optimization problems is to base the
preconditioner on the solution from a reduced approximated problem
where the Schur complement is replaced by an approximation \eg{} by
quasi-newton methods, see \cite{ghattas}.

In our case we use the GMRES method 
to solve the non-symmetric system \eqref{eq:saddle_syst} and base our
preconditioner on the approximate solution of a simple blockwise
Gauss-Seidel method
\ie{} to start with a guess $\hat q^0$ and iteratively solve 
\begin{equation}
  \label{eq:GS}
  \begin{aligned}
    K_{11} \hat u^{i+1} &= f - K_{12} \hat q^i,\\
    K_{11}^T \hat q^{i+1} &= g - K_{21} \hat u^{i+1},
  \end{aligned}
\end{equation}
which works well for large regularizations \ie{} when $\mathfrak{h}_\delta''$
is small and the diagonal blocks of \eqref{eq:saddle_syst} are dominant. 
Also, each iteration with this method only requires one forward and
one backward solve in time of a modified heat equation so the
computational work for one iteration is concentrated to solving $N-1$
smaller systems with system matrices $(M+S_i)$. 
In practice, the Gauss-Seidel method will break down for small
regularizations but for our problems (and discretizations) only one
iteration with \eqref{eq:GS} turns out to be a fairly good
approximation to use as preconditioner. Note that for $\hat
q^0=0$, one Gauss-Seidel iteration is the same as solving
\eqref{eq:saddle_syst} with the approximation $K_{12}=0$.

% To solve the original system \eqref{eq:saddle_syst} for
% $\mathfrak{h}''=0$ means to solve the heat
% equation forward and backward in time, \ie{} $K_{11}$ can interpreted as
% the differential operator $\partial_t-\Delta$, and to find a good
% approximation to the Schur complement 
% \begin{equation*}
%   K_{11}^T -  K_{21} K_{11}^{-1} K_{12}, 
% \end{equation*}
% for the \eqref{eq:saddle_syst} is essentially to find
% a good approximation for the inverse of $K_{11}$.

Another more elaborate idea is to use a preconditioner based on
the solution of an approximated Schur complement
system 
\begin{equation*}
  \begin{aligned}
   \left(
    \begin{array}{cc}
      K_{11} & K_{12}\\
      0 & S
    \end{array}
  \right)
  \left(
    \begin{array}{c}
      \hat u\\
      \hat q
    \end{array}
  \right)
  =
  \left(
    \begin{array}{c}
      g\\
      f - K_{12} K_{11}^{-1} g
    \end{array}
  \right),
   \end{aligned}
 \end{equation*}
where $S$ is an approximation of the Schur complement
\begin{equation*}
  K_{11}^T - K_{21} K_{11}^{-1} K_{12}.
\end{equation*}
which essentially is to find a good approximation of the lower
triangular block matrix $K_{11}^{-1}$.

Although solution algorithms for saddle point systems on the symmetric
form \eqref{eq:saddle_syst2} are
extensively treated in the litterature, see \cite{saddle-point} for
an overview, we here favour the 
non-symmetric form \eqref{eq:saddle_syst}, since a Schur complement
reduction of \eqref{eq:saddle_syst2} means to find an approximation to
the Schur complement
\begin{equation*}
  K_{12} - K_{11} K_{21}^{-1} K_{11}^T,
\end{equation*}
which since $K_{21}$ here can be singular, is unavailable. One way
around this obstacle is to rewrite \eqref{eq:saddle_syst2} by \eg{}
the augmented Lagrangian method which leads to a symmetric invertible 
Schur complement but where the physical meaning of the original
system, on PDE level, is partially lost.

If a direct solver is used for the Newton system it is appropriate to 
reorder \eqref{eq:saddle_syst} such that the solution vector and right
hand side contains time steps in increasing order, which leads to a
banded Jacobian with band-width of the same order as the number of
spacial degrees of freedom. 

% One way around this obstacle is to rewrite \eqref{eq:saddle_syst2} 
% by \eg{} the augmented Lagrangian method which introduces a
% symmetric positive semidefinite matrix $W$, \eg{} $W=\gamma I$ for some
% $\gamma>0$, and the equivalent system
% \begin{equation*}
%     \left(
%     \begin{array}{cc}
%       K_{21} + K_{11}^T W K_{11} & K_{11}^T + K_{11}^T W K_{12}\\
%       K_{11} & K_{12}
%     \end{array}
%   \right)
%   \left(
%     \begin{array}{c}
%       \hat u\\
%       \hat q
%     \end{array}
%   \right)
%   =
%   \left(
%     \begin{array}{c}
%       g +  K_{11} W f\\
%       f
%     \end{array}
%   \right),
% \end{equation*}
% which unfortunately is not symmetric but where $K_{21} + K_{11}^T W
% K_{11}$ is invertible and positive definite for a suitable
% choice of $W$. To construct a preconditioner, the Schur complement
% \begin{equation*}
%   K_{12} -  K_{11} (K_{21} + K_{11}^T W K_{11})^{-1} (K_{11}^T
%   + K_{11}^T W K_{12} ), 
% \end{equation*}
% needs to be replaced by an efficiently calculated approximation, which
% is a delicate task since the physical meaning of the original system,
% on PDE level, is partially lost. 

Our computations were implemented MATLAB (for the one dimensional
examples), and in DOLFIN \cite{dolfin}, the C++/Python interface of
the finite element solver environment FEniCS \cite{fenics} (for the
two dimensional examples). Piecewise linear basis functions were used
for the finite element subspace $\bar V$, and in all examples the 
solution $u,q$ was first calculated for a large regularization which
was succesively reduced such that the solution from the previous
regularization served as starting guess for a smaller regularization.

For the two dimensional
examples the sadde-point system
\eqref{eq:saddle_syst} was solved with 
the PETSc implementation of GMRES (used by DOLFIN)
with preconditioning based on the solution from one iteration of 
blockwise Gauss-Seidel method. For the one dimensional examples a
direct solver was used. The number of iterations for GMRES
with the Gauss-Seidel preconditioner seems to be relatively insensitive with
respect to temporal and spacial discretization but still highly sensitive to
the regularization in our examples.

% TODO: Explain the implications for solving saddle-point systems
% efficiently. Explain averaging and why time independent control cannot be used.

To give a time independent approximation $\sigma(x)$ of the time dependent control
$\sigma(x,t)$, approximated by $\tilde \sigma:=\mathfrak{h}_\delta'(\nabla u \cdot
\nabla q)$ where $u,q$ are solutions to the Hamiltonian system
\eqref{eq:hamilt_syst_heat}, three different types of averaging were
tested as post-processing: 
\begin{enumerate}

\item 
Let the time independent control be defined by the Hamiltonian 
\eqref{eq:time_indep_hamilt}, \ie{}
\begin{equation}
  \label{eq:avg1}
  \sigma:=\mathfrak{h}_\delta'\bigg( \int_0^T \nabla u \cdot \nabla q \ud
  t \bigg).
\end{equation}

\item
Let the time independent control be the average of the time dependent
control, \ie{}
\begin{equation}
  \label{eq:avg2}
  \sigma:=\frac{1}{T} \int_0^T \mathfrak{h}_\delta'(\nabla u \cdot
  \nabla q) \ud t.
\end{equation}

\item
Let the time independent control be the weighted average 
\begin{equation}
  \label{eq:avg3}
  \sigma:=\frac{ \int_0^T \mathfrak{h}_\delta'(\nabla u \cdot
    \nabla q) |\nabla u \cdot \nabla q| \ud t }
  { \int_0^T  |\nabla u \cdot \nabla q| \ud t },
\end{equation}
of the time dependent control $\mathfrak{h}_\delta'(\nabla u \cdot\nabla q)$.

\end{enumerate}
The weighted average turned out to be the most successful aproximation
and can be explained by first extending the Hamiltonian
\eqref{eq:hamilt_heat} to also depend 
on the artifical variable $z$ as in Section
\ref{sec:time-indep-control}
\begin{equation*}
  \begin{aligned}
    H(u,q,z) :=& 
    \frac{1}{\tilde T} \int_0^T \int_{\partial\Omega} (u-u^*)^2 + j q \ud s \ud t - 
    \frac{1}{\tilde T} \int_0^T \int_\Omega u_t q \ud x \ud t\\
    & -\frac{1}{\tilde T} \int_\Omega \int_0^T 
     \mathfrak{h}'(\nabla u \cdot \nabla q) \nabla u \cdot \nabla q \ud t \ud x,\\
  \end{aligned}
\end{equation*}
where $\mathfrak{h}'(\nabla u \cdot \nabla q) \nabla u
\cdot \nabla q = \mathfrak{h}(\nabla u \cdot \nabla q)$ by definition.
For the problem with a time independent control we now seek an
approximation of the Hamiltonian \eqref{eq:time_indep_hamilt} of the form
\begin{equation*}
  \begin{aligned}
     \bar H(u,q,z) :=& 
     \frac{1}{\tilde T} \int_0^T \int_{\partial\Omega} (u-u^*)^2 + j q \ud s \ud t - 
    \frac{1}{\tilde T} \int_0^T \int_\Omega u_t q \ud x \ud t\\
    & -\frac{1}{\tilde T} \int_\Omega
     f( \nabla u \cdot \nabla q)
     \int_0^T \nabla u \cdot \nabla q  \ud t \ud x,\\
  \end{aligned}
\end{equation*}
that best approximates $H$, \ie{}
\begin{equation*}
    f(\nabla u \cdot \nabla q):=\frac{ \int_0^T \mathfrak{h}'(\nabla u \cdot
      \nabla q) \nabla u \cdot \nabla q \ud t }
    { \int_0^T  \nabla u \cdot \nabla q \ud t }.
\end{equation*}

\begin{figure}[htbp]
  \centering
  % \psfrag{x}{$x$}
  % \psfrag{t}{$t$}
  % \includegraphics[height=0.25\textwidth]{./pictures/heat/test1_same_mesh/u_plot.eps}  
  % \includegraphics[height=0.25\textwidth]{./pictures/heat/test1/u_plot.eps}  
  % \includegraphics[height=0.25\textwidth]{./pictures/heat/test1_noise/u_plot.eps}\\ 
  % \includegraphics[height=0.25\textwidth]{./pictures/heat/test1_same_mesh/p_plot.eps}  
  % \includegraphics[height=0.25\textwidth]{./pictures/heat/test1/p_plot.eps}  
  % \includegraphics[height=0.25\textwidth]{./pictures/heat/test1_noise/p_plot.eps}\\ 
  % \includegraphics[height=0.25\textwidth]{./pictures/heat/test1_same_mesh/c_rec.eps}
  % \includegraphics[height=0.25\textwidth]{./pictures/heat/test1/c_rec.eps}
  % \includegraphics[height=0.25\textwidth]{./pictures/heat/test1_noise/c_rec.eps}\\
  % \psfrag{x}{$\delta$}
  % \includegraphics[height=0.25\textwidth]{./pictures/heat/test1_same_mesh/obj.eps}  
  % \includegraphics[height=0.25\textwidth]{./pictures/heat/test1/obj.eps}  
  % \includegraphics[height=0.25\textwidth]{./pictures/heat/test1_noise/obj.eps
  \includegraphics[width=1\textwidth]{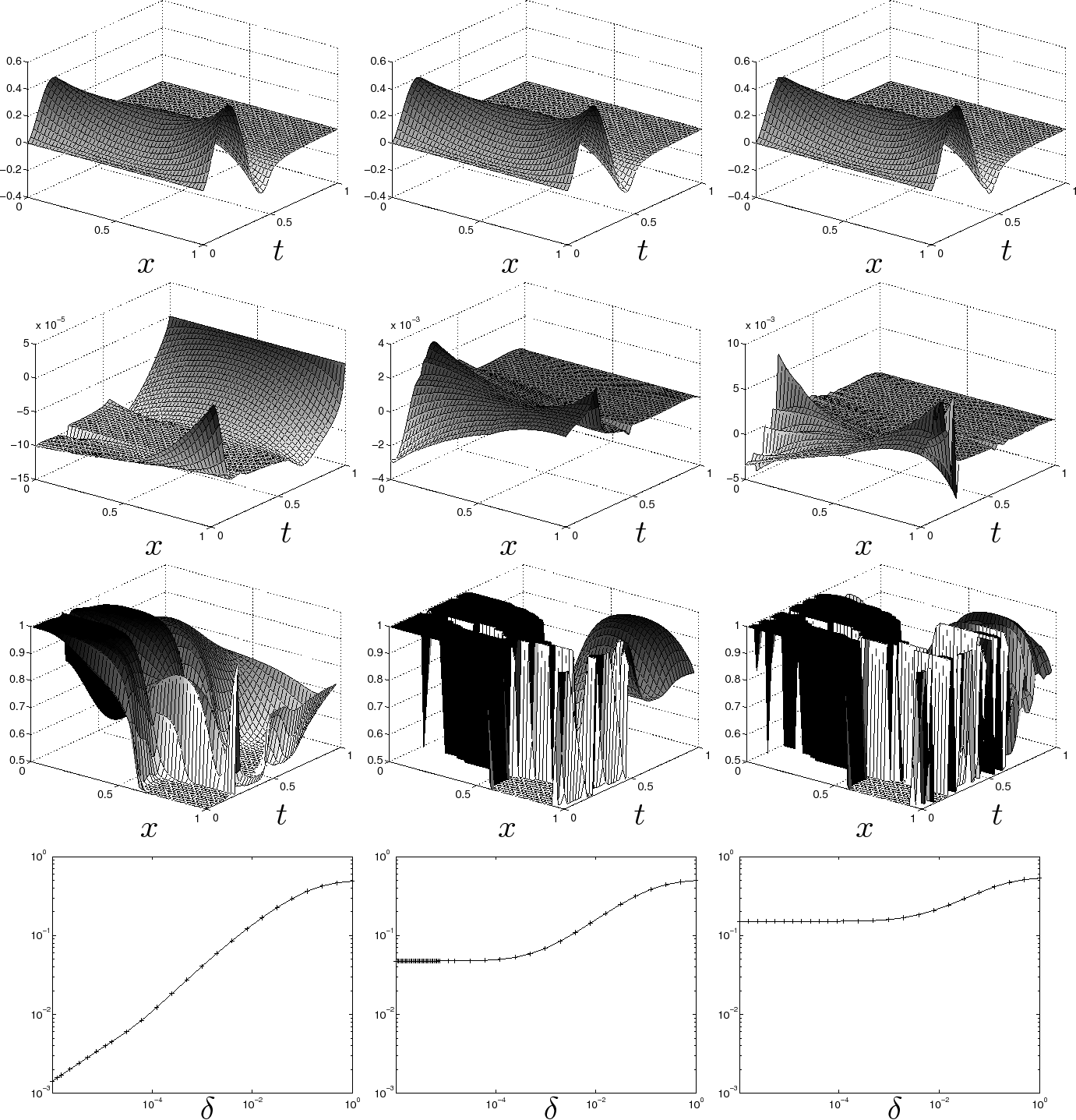}
  \caption{1D reconstruction of $\sigma_{true} =
    0.75-0.5\tanh(20x-10)$ with $\delta=10^{-6}$, $\sigma_-=0.5$,
    $\sigma_+=1$, 
    measurements on both boundaries and Neumann boundary condition 
    $\sigma u_x(0,t) = -\sigma u_x(1,t) = \sin(4t)$ for $t<0.5$ and
    $0$ elsewhere.  
    The plot shows, from top to bottom, $u$, $q$,
    $\mathfrak{h}_\delta'$ and the objective function
    $\|u-u^*\|_{L^2(\partial\Omega\times [0,T])}$. 
    In all cases $u, q$ was
    calculated with 50 steps in space and time.
    In the left column, the data $u^*$ was generated by solving the
    heat equation with 50 steps in time and space and conductivity
    $\sigma_{true}$, while 200 steps in time and space was used in the
    middle and right columns.
    In the right column 10\% noise was also added to $u^*$.
  }
  \label{fig:constr1}
\end{figure}

In Figure \ref{fig:constr1}, one dimensional
reconstructions from three sets of simulated data $u^*$, generated from
a time independent conductivity $\sigma_{true}$, are compared:  
\begin{enumerate}
\item Data calculated with the same discretization as $u$ and $q$.
\item Different discretisations used for data and solutions.
\item Different discretisations used for data and solutions and
  with noise in the data $u^*$. 
\end{enumerate}
The last set is the most realistic one since for true experimental
data of $u^*$ it is inevitable to not only have measurement noise but
also a systematic error from the numerical method. 
To simulate noise the discrete
solution $u^*$ was multiplied componentwise by 
independent standard normal distributed stochastic variables
$\varepsilon_{ij}$ according to $u^*(x_i,t_j)(1+\eta \epsilon_{ij})$, 
where $\eta$ denotes the percentage of noise. 
It is notable that the
systematic error from using different meshes  
can have a much bigger effect on the solutions than additional
noise, which can be observed from the dual solution $q$ in
Figure \ref{fig:constr1}.  

\begin{figure}[htbp]
  \centering
  \includegraphics[width=1\textwidth]{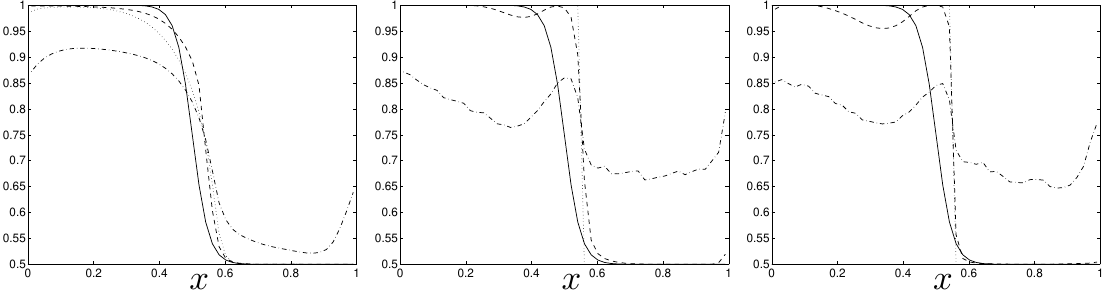}
  \caption{The time independent post-processed conductivity for the 1D
    reconstructions in Figure \ref{fig:constr1}. The true control
    $\sigma_{true}$ is indicated by a solid line and the averaged
    controls \eqref{eq:avg1}, \eqref{eq:avg2} and \eqref{eq:avg3} are
    indicated by dotted, dash-dotted and dashed lines, respectively.  
  }
  \label{fig:cond1}
\end{figure}

\begin{figure}[htbp]
  \centering
  % \psfrag{x}{$x$}
  % \psfrag{t}{$t$}
  % \includegraphics[height=0.25\textwidth]{./pictures/heat/test2_same_mesh/c_rec.eps}
  % \includegraphics[height=0.25\textwidth]{./pictures/heat/test2/c_rec.eps}
  % \includegraphics[height=0.25\textwidth]{./pictures/heat/test2_noise/c_rec.eps}\\
  % \includegraphics[height=0.25\textwidth]{./pictures/heat/test2_same_mesh/c_rec_avg.eps}
  % \includegraphics[height=0.25\textwidth]{./pictures/heat/test2/c_rec_avg.eps}
  % \includegraphics[height=0.25\textwidth]{./pictures/heat/test2_noise/c_rec_avg.eps}\\
  \includegraphics[width=1\textwidth]{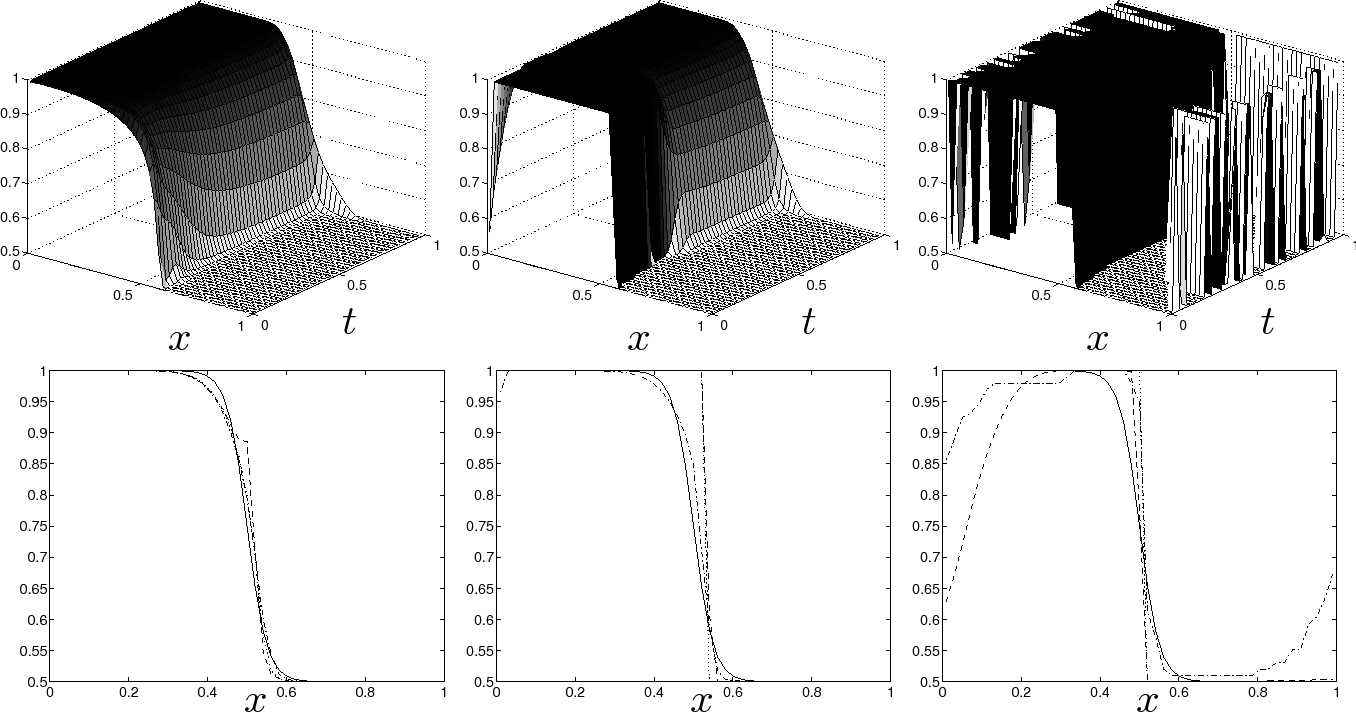}
  \caption{
    1D reconstruction with data as in Figure \ref{fig:constr1}
    and \ref{fig:cond1} but with
    Neumann boundary condition 
    $\sigma u_x(0,t) = -\sigma u_x(1,t) = 1$.
  }
  \label{fig:cond2}
\end{figure}

In Figure \ref{fig:cond1} the time independent post-processing of the time
dependent reconstruction can be found.
It is here evident that the weighted average \eqref{eq:avg3} performs
better than \eqref{eq:avg2}, but since the reconstruction is highly dependent on the
given boundary condition, see Figure \ref{fig:cond2} for comparison,
there are situations where the different post-processing techniques
perform equally well. It would of course be optimal to use the
knowledge that $\sigma_{true}$ is independent of time in the
calculations, \ie{} to use the Newtion system for
\eqref{eq:hamilt_syst_heat} 
with time independent-control \eqref{eq:control_time_indep}. This
would however lead to a dense Jacobian.  

Note that in the examples the limits
$\sigma_-,\sigma_+$ were chosen to be the biggest and smallest values
of $\sigma_{true}$. In our experience the Pontryagin method is not
well suited for reconstruction of values between
$\sigma_-$ and $\sigma_+$ if there is noise or other measurement errors
present in data.

\begin{figure}[htbp]
  \centering
  % \psfrag{x}{$x$}
  % \psfrag{y}{$y$}
  % \includegraphics[height=0.25\textwidth]{pictures/heat2d/test1/trueCoeff.eps}
  % \includegraphics[height=0.25\textwidth]{pictures/heat2d/test1/avgCoeff.eps}
  % \includegraphics[height=0.25\textwidth]{pictures/heat2d/test1_noise/avgCoeff.eps}\\
  % \includegraphics[height=0.25\textwidth]{pictures/heat2d/test2/trueCoeff.eps}
  % \includegraphics[height=0.25\textwidth]{pictures/heat2d/test2/avgCoeff.eps}
  % \includegraphics[height=0.25\textwidth]{pictures/heat2d/test2_noise/avgCoeff.eps}
  \includegraphics[width=1\textwidth]{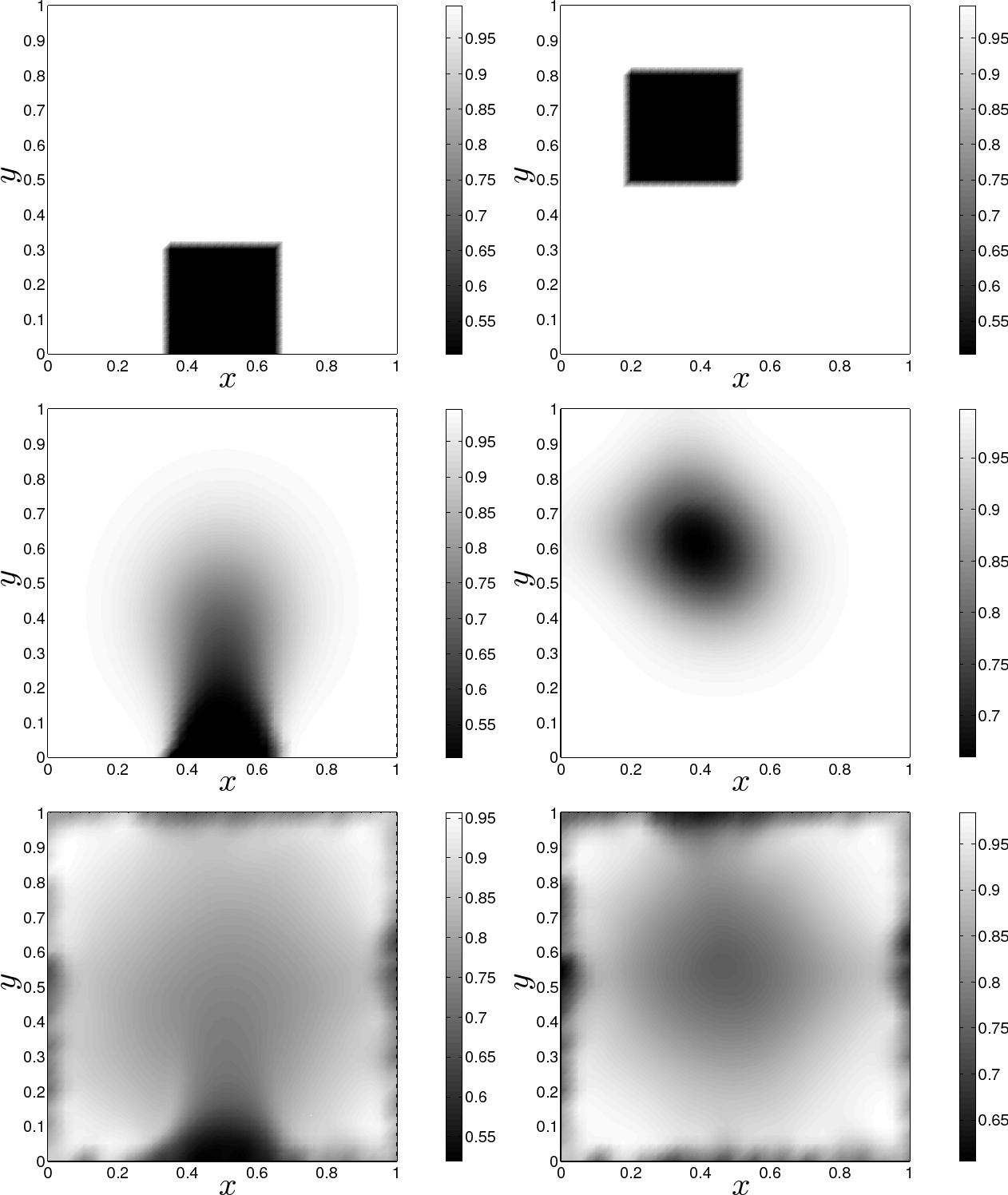}
  \caption{2D reconstruction on the unit square with final time $T=1$ and
    Neumann boundary condition 
    $\sigma \frac{\partial u}{\partial n} = 1$ on $\partial \Omega
    \times [0,T]$.
    The data $u^*$ was simulated by solving the forward equation on a
    quasi-uniform mesh with 13000 triangles and 80 time steps while the
    inverse problem was solved on a uniform mesh with 3200
    triangles and 40 time steps. 
    Measurements from the whole boundary were used.
    Top: True conductivity $\sigma_{true}$. Middle: Reconstructed
    condictivity for $\delta\approx 0.002$ using the weighted average
    \eqref{eq:avg3}. 
    Bottom: As in middle but for $\delta\approx 0.05$ and with 5\%
    noise in the measurements.
  }
  \label{fig:2Dheat}
\end{figure}

Figure \ref{fig:2Dheat} shows two-dimensional reconstructions of
two different time independent conductivities. Unlike the
one-dimensional example the quality of the reconstruction here deteriorates
quickly as the distance to the measurement locations is increased.

\section{Reconstruction from the Wave Equation }\label{sec:wave}

In this section the goal is to determine the wave speed for a scalar
acoustic wave equation: 
Given measured data $u^*$, find the state $u : \bar \Omega \times [0,T] \to
V$, $u=u(x,t)$ and
a control \eg{} $\sigma : \bar \Omega \times [0,T] \to
[\sigma_-,\sigma_+]$, $\sigma=\sigma(x,t)$ where $0<\sigma_-<\sigma_+$, that solves the
partial differential equation 
\begin{equation} 
  \label{eq:constr_wave}
  \begin{aligned} 
    u_{tt} &= \divop ( \sigma \nabla u ),
    & \text{ in } \Omega \times (0,T],\\
    \sigma \nabla u \cdot \bf{n} &= j, 
    & \text{ on } \partial\Omega \times (0,T],\\
    u &= u_t = 0,
    & \text{ on } \bar \Omega \times \{t=0\},
  \end{aligned}
\end{equation}
such that the error functional
\begin{equation} 
  \label{eq:obj_wave}
  \int_0^T \int_{\partial\Omega} (u-u^*)^2 \ud s \ud t,
\end{equation}
is minimized. 
The control $\sigma$ is here the square of the wave speed of the
medium and $u$ is the pressure deviation. 

To use the framework of the previous section we note that
\eqref{eq:constr_wave} can be written as the first order system
\begin{equation}
  \label{eq:constr_wave_syst}
  \begin{aligned}
    v_{t} &= \divop ( \sigma \nabla u ), & \text{ in } \Omega \times (0,T],\\
    u_t &= v, & \text{ in } \bar \Omega \times (0,T],\\
    \sigma \nabla u \cdot \bf{n} &= j, 
    & \text{ on } \partial\Omega \times (0,T],\\
    u &= v = 0,
    & \text{ on } \bar \Omega \times \{t=0\}.
  \end{aligned}
\end{equation}

\subsection{The Hamiltonian System}
As in Section \ref{sec:hamilt_heat} we have a Hamiltonian
associated with the optimal 
control problem \eqref{eq:obj_wave} and \eqref{eq:constr_wave_syst}
which is defined by
\begin{equation}
  \label{eq:hamilt_wave}
  \begin{aligned}
  {H} &:= \min_{\sigma : \Omega \to [\sigma_-,\sigma_+]} 
  \int_{\partial\Omega} (u-u^*)^2 \ud s + 
  \int_\Omega \divop ( \sigma \nabla u ) q + vp \ud x\\
  &=\int_{\partial\Omega} (u-u^*)^2 + j q \ud s + 
  \min_{\sigma : \Omega \to [\sigma_-,\sigma_+]} 
   \int_\Omega - \sigma \nabla u \cdot \nabla q + vp \ud x\\
   &=\int_{\partial\Omega} (u-u^*)^2 + j q \ud s + 
   \int_\Omega vp - 
   \underbrace{
     \max_{\sigma \in [\sigma_-,\sigma_+]} \{\sigma \nabla u \cdot \nabla q\} 
   }_{\mathfrak{h}(\nabla u \cdot \nabla q)} 
   \ud x,
  \end{aligned}
\end{equation}
and the Hamiltonian system becomes
\begin{equation}
  \label{eq:hamilt_syst_wave}
  \begin{aligned}
    v_t &= \divop \big( \tilde \sigma \nabla
    u \big), & \text{ in } \Omega \times (0,T],\\
    u_t &= v, & \text{ in } \bar \Omega \times (0,T],\\
    \tilde \sigma \nabla u \cdot \bf{n} &= j, 
    & \text{ on } \partial\Omega \times (0,T],\\
    u &= v = 0, & \text{ on } \bar \Omega \times \{t=0\},\\
    -p_t &= \divop \big( \tilde \sigma \nabla
    q \big), & \text{ in } \Omega \times (0,T],\\
    -q_t &= p, & \text{ in } \bar \Omega \times (0,T],\\
    \tilde \sigma \nabla q \cdot \bf{n} &=  2(u-u^*), 
    & \text{ on } \partial\Omega \times (0,T],\\
    p &= q = 0, & \text{ on }\Omega \times \{t=T\},\\
  \end{aligned}
\end{equation}
or equivalently
\begin{equation}
  \label{eq:hamilt_syst_wave2}
  \begin{aligned}
    u_{tt} &= \divop \big( \tilde \sigma \nabla u \big),
    & \text{ in } \Omega \times (0,T],\\
    \tilde \sigma 
    \nabla u \cdot \bf{n} &= j, & \text{ on } \partial\Omega \times (0,T],\\
    u &= u_t = 0,
    & \text{ on } \bar\Omega \times \{t=0\},\\
    q_{tt} &= \divop \big( \tilde \sigma \nabla q \big),
    & \text{ in } \Omega \times (0,T],\\
    \tilde \sigma 
    \nabla q \cdot \bf{n} &= 2(u-u^*), 
    & \text{ on } \partial\Omega \times (0,T],\\
    q &= q_t = 0,
    & \text{ on } \bar \Omega \times \{t=T\}.\\
  \end{aligned}
\end{equation}
with 
\begin{equation*}
  \tilde \sigma := \mathfrak{h}'(\nabla u \cdot \nabla q).
\end{equation*}

\subsection{Symplecticity for the Wave Equation}
As a natural case the symplectic methods discussed in
\ref{sec:symplectic_heat}, with $\varphi = (u, v)$, $\lambda = (p, q)$, can be
used to solve the system \eqref{eq:hamilt_syst_wave}. It is however
also possible to
use a time-discretization that is symmetric in time i.e.
\begin{equation}
  \label{eq:sympl_wave}
  \begin{aligned}
    u_{n+1}-2u_n+u_{n-1} &= k^2 \divop 
    \big( \tilde \sigma_n \nabla u_{n}
    \big), & \text{ in } \Omega,\\
    \tilde \sigma_n \nabla u_{n} \cdot
    \bf{n} &= j_{n}, & \text{ on } \partial \Omega,\\
    u_0 &= u_1 = 0, & \text{ in } \Omega,\\
    q_{n+1}-2q_n+q_{n-1} &= k \divop 
    \big( \tilde \sigma_n \nabla
    q_{n} \big), & \text{ in } \Omega,\\
    \tilde \sigma_n \nabla q_{n} \cdot
    \bf{n} &= 2(u_{n}-u^*_{n}), & \text{ on } \partial \Omega,\\
    q_N & = q_{N-1} = 0, & \text{ in } \Omega,\\
  \end{aligned}
\end{equation}
for $\tilde \sigma_n := \mathfrak{h}'(\nabla u_n \cdot \nabla q_n)$ and $n=1,\ldots,N-1$.
For a given $\tilde \sigma$, constant in time, this scheme is
the symplectic backward Euler
method for the forward wave equation for $u$, which can be written as the
Hamiltonian system \eqref{eq:constr_wave_syst} with Hamiltonian
\begin{equation*}
  {H}_{wave}(u,v) := \frac{1}{2} \int_\Omega |\tilde \sigma \nabla u|^2 +
  v^2 \ud x,
\end{equation*}
and the symplectic forward Euler
method for the backward wave equation for $q$.

To see that that the scheme \eqref{eq:sympl_wave} is symplectic for 
$\tilde \sigma_n := \mathfrak{h}'(\nabla u_n \cdot \nabla q_n)$ we note that a
one-step method $(\varphi_n,\lambda_n) \to 
(\varphi_{n+1},\lambda_{n+1})$ is symplectic if there exists a
function ${H}(\varphi_n,\lambda_{n+1})$ such that
\eqref{eq:sympl_forw_euler} holds, or equivalently
${H}(\varphi_{n+1},\lambda_{n})$ such that
\eqref{eq:sympl_backw_euler} holds, see Remark 4.8 in \cite{ss} or
\cite{wanner} for details. It thus follows that the one-step method
\begin{equation*}
  \begin{aligned}
    v_{n+1}-v_{n} &= k \divop 
    \big(  \mathfrak{h}'(\nabla u_{n} \cdot \nabla q_{n}) \nabla u_{n}
    \big),\\
    u_{n+1}-u_{n} &= k v_{n+1},\\
    p_{n}-p_{n+1} &= k \divop 
    \big( \mathfrak{h}'(\nabla u_{n} \cdot \nabla q_{n}) \nabla
    q_{n} \big),\\
    q_{n}-q_{n+1} &= k p_{n+1},\\
  \end{aligned}
\end{equation*}
corresponds to the symplectic forward Euler method for the Hamiltonian
\begin{equation*}
  \begin{aligned}
    \tilde H (
    \underbrace{u_{n},q_{n}}_{\varphi_n},
    \underbrace{v_{n+1},p_{n+1}}_{\lambda_{n+1}}) := 
    H(u_{n},v_{n+1},p_{n+1},q_{n}) - 2\int_\Omega v_{n+1}p_{n+1} \ud x,
  \end{aligned}
\end{equation*}
where $H$ is given by \eqref{eq:hamilt_wave}. 
Since \eqref{eq:sympl_wave} only is stable for sufficiently
small time-steps and still requires to solve a complex saddle point system we
will use the symplectic midpoint method in
our experiments.

\subsection{Numerical Examples}
Let $\tilde \sigma:= \mathfrak{h}'_\delta$
where $\mathfrak{h}'_\delta$ is given by \eqref{eq:reg}.
The symplectic midpoint method for the regularized Hamiltonian system
\eqref{eq:hamilt_syst_wave} can then be written as 
\begin{equation*}
  \begin{aligned}
    F^1_n(w) = 0, \quad 
    F^2_n(w) = 0, \quad      
    G^1_n(w) = 0, \quad 
    G^2_n(w) = 0, \quad  
  \end{aligned}
\end{equation*}
for $n=0,\ldots,N-1$, and $ \forall w \in \bar V$,
where
\begin{equation*}
  \begin{aligned}
    F^1_n(w) &:= \int_\Omega (v_{n+1}-v_{n})w + k \mathfrak{h}'_\delta\big(
    \nabla u_{n+\frac{1}{2}} \cdot \nabla q_{n+\frac{1}{2}}
    \big) \nabla u_{n+\frac{1}{2}} \cdot \nabla w \ud x\\
    &- \int_{\partial \Omega} k j_{n+\frac{1}{2}} w \ud s,\\
    F^2_n(w) &:= \int_\Omega ( u_{n+1}-u_{n} - k v_{n+\frac{1}{2}}) w \ud x,\\
    G^1_n(w) &:= \int_\Omega (q_{n}-q_{n+1} - k p_{n+\frac{1}{2}})w \ud x.\\
    G^2_n(w) &:= \int_\Omega (p_{n}-p_{n+1})w + k \mathfrak{h}'_\delta(\nabla 
    u_{n+\frac{1}{2}} \cdot \nabla q_{n+\frac{1}{2}}) \nabla q_{n+\frac{1}{2}} \cdot \nabla w \ud x\\
    &-  \int_{\partial \Omega} 2k (u_{n+\frac{1}{2}}-u^*_{n+\frac{1}{2}}) w \ud s,\\
  \end{aligned}
\end{equation*}
and $u_0=v_0=p_N=q_N=0$.  The index $n+\frac{1}{2}$ implies the
average of the values at $n$ and $n+1$,
\ie{} $u_{n+\frac{1}{2}}:=\frac{1}{2}(u_n+u_{n+1})$.  Taking the
variations with respect to 
$u,v,p,q$ gives the Newton system 
\begin{equation}
  \label{eq:saddle_syst_wave}
    \left(
    \begin{array}{cccc}
      K_{11} & K_{12} & 0 & K_{14}\\
      K_{21} & K_{22} & 0 & 0\\
      0 & 0 & K_{33} & K_{34}\\
      K_{41} & 0 & K_{43} & K_{44}\\
    \end{array}
  \right)
  \left(
    \begin{array}{c}
      \hat u\\
      \hat v\\
      \hat p\\
      \hat q
    \end{array}
  \right)
  = 
  \left(
    \begin{array}{c}
      f_1\\
      f_2\\
      g_1\\
      g_2\\
    \end{array}
  \right),
\end{equation}
with increments
\begin{equation*}
  \begin{aligned}
    \hat u = &
    \left(
      \begin{array}{ccccc}
        \hat u_1 & \ldots & \hat u_N
      \end{array}
    \right)^T, 
    & \hat v = &
    \left(
      \begin{array}{ccccc}
        \hat v_1 & \ldots & \hat v_{N}
      \end{array}
    \right)^T,\\
    \hat p = &
    \left(
      \begin{array}{ccccc}
        \hat p_0 & \ldots & \hat p_{N-1}
      \end{array}
    \right)^T, 
    & \hat q = &
    \left(
      \begin{array}{ccccc}
        \hat q_0 & \ldots & \hat q_{N-1}
      \end{array}
    \right)^T,\\
  \end{aligned}
\end{equation*}
and right hand side
\begin{equation*}
  \begin{aligned}
    f_1 = &
    \left(
      \begin{array}{ccccc}
        F^1_0 & \ldots & F^1_{N-1}
      \end{array}
    \right)^T,
    & f_2 = &
    \left(
      \begin{array}{ccccc}
        F^2_0 & \ldots & F^2_{N-1}
      \end{array}
    \right)^T,\\
    g_1 = &
    \left(
      \begin{array}{ccccc}
        G^1_0 & \ldots & G^1_{N-1}
      \end{array}
    \right)^T,
    & g_2 = &
    \left(
      \begin{array}{ccccc}
        G^2_0 & \ldots & G^2_{N-1}
      \end{array}
    \right)^T.\\
\end{aligned}
\end{equation*}
with submatrices with the following structure:
\begin{itemize}
\item $K_{11}$ is lower block bi-diagonal with 
  \begin{equation}
    \label{eq:K11}
    \begin{aligned}
      & \frac{1}{2} \int_\Omega k \mathfrak{h}''_\delta(\nabla
      u_{n+\frac{1}{2}} \cdot \nabla q_{n+\frac{1}{2}}) \nabla q_{n+\frac{1}{2}} \cdot \nabla w \
      \nabla u_{n+\frac{1}{2}} \cdot \nabla \bar w \ud x\\
      &+\frac{1}{2} \int_\Omega k \mathfrak{h}'_\delta(\nabla
      u_{n+\frac{1}{2}} \cdot \nabla q_{n+\frac{1}{2}}) \nabla w \cdot \nabla \bar w \ud x,\\
    \end{aligned}
  \end{equation}
  on its main diagonal for $n=0,\ldots,N-1$ and on its sub-diagonal
  for $n=1,\ldots,N-1$.

\item $K_{44}$ is upper block bi-diagonal with \eqref{eq:K11} on its
  diagonal for $n=0,\ldots,N-1$ and on its super-diagonal for $n=0,\ldots,N-2$.
  
\item $K_{12}=K_{21}=K_{34}^T=K_{43}^T$ is lower block bi-diagonal with mass matrices
  $M$ on the main diagonal and $-M$ on the subdiagonal. 

\item $K_{22}=K_{33}^T$ is lower block bi-diagonal with
  $-\frac{kM}{2}$ on the diagonal and the sub-diagonal. 

\item $K_{14}$ is upper block bi-diagonal with
  \begin{equation*}
    \frac{1}{2} \int_\Omega k \mathfrak{h}''_\delta(\nabla
    u_{n+\frac{1}{2}} \cdot \nabla q_{n+\frac{1}{2}}) \nabla u_{n+\frac{1}{2}} \cdot \nabla w \
    \nabla u_{n+\frac{1}{2}} \cdot \nabla \bar w \ud x,
  \end{equation*}
  on its diagonal for $n=0,\ldots,N-1$ and on its super-diagonal for $n=0,\ldots,N-2$. 

\item $K_{41}$ is lower block bi-diagonal with
  \begin{equation*}
    \frac{1}{2} \int_\Omega k \mathfrak{h}''_\delta(\nabla
    u_{n+\frac{1}{2}} \cdot \nabla q_{n+\frac{1}{2}}) \nabla q_{n+\frac{1}{2}} \cdot \nabla w \
    \nabla q_{n+\frac{1}{2}} \cdot \nabla \bar w \ud x -
    \int_{\partial\Omega} k \bar w w \ud s,
  \end{equation*}
  on its diagonal for $n=0,\ldots,N-1$ and sub-diagonal for $n=1,\ldots,N-1$. 

\end{itemize}

As in the previous section we will solve the Newton system using GMRES
and an approximate solution as preconditioner, \eg{} from the 
% $4 \times 4$ blockwise Gauss-Seidel method
% \begin{equation}
%   \label{eq:GSwave}
%   \begin{aligned}
%     K_{11} \hat u^{n+1} &= f_1 - K_{12} \hat v^n - K_{14} \hat q^n,\\
%     K_{22} \hat v^{n+1} &= f_2 - K_{21} \hat u^{n+1},\\
%     K_{33} \hat p^{n+1} &= g_1 - K_{34} \hat q^{n},\\
%     K_{44} \hat q^{n+1} &= g_2 - K_{41} \hat u^{n+1} - K_{43} \hat p^{n+1},\\
%   \end{aligned}
% \end{equation}
% or
the $2 \times 2$ blockwise Gauss-Seidel method
\begin{equation*}
  \begin{aligned}
    K_{11} \hat u^{i+1} + K_{12} \hat v^{i+1} &= f_1  - K_{14} \hat q^i,\\
    K_{21} \hat u^{i+1} + K_{22} \hat v^{i+1} &= f_2,\\
    K_{33} \hat p^{i+1} + K_{34} \hat q^{i+1} &= g_1,\\
    K_{43} \hat p^{i+1} + K_{44} \hat q^{i+1} &= g_2 - K_{41} \hat u^{i+1},\\
  \end{aligned}
\end{equation*}
which can be written as
\begin{equation}
  \label{eq:GSwave2}
  \begin{aligned}
    ( K_{11} - K_{12}K_{22}^{-1}K_{21} ) \hat u^{i+1} &= f_1 -
    K_{12}K_{22}^{-1}f_2 - K_{14} \hat q^i,\\ 
    ( K_{44} - K_{43}K_{33}^{-1}K_{34} ) \hat q^{i+1} &= g_2 -
    K_{43}K_{33}^{-1} g_1 - K_{41} \hat u^{i+1}.\\
  \end{aligned}
\end{equation}
Note that \eqref{eq:GSwave2} is easily
solved since inverting
$K_{22}$ and $K_{33}$ only involves the calculation of $M^{-1}$. In fact,
the Schur complements $K_{11} - K_{12}K_{22}^{-1}K_{21}$ and
$K_{44}-K_{43}K_{33}^{-1}K_{34}$ becomes 
lower and upper block trianglar matrices, respectively, and \eqref{eq:GSwave2}
can be solved by
one forward substitution in time for $\hat u^{i+1}$ and one backward
substitution in time for $\hat q^{i+1}$. Of course, to save memory the
Schur complement system \eqref{eq:GSwave2} should never be formed explicitly.
For large regularizations the
Schur complements can be seen as approximations of the operator
$-\Delta + \partial_{tt}$.
As for the case with the heat
equation starting with $\hat q^0=0$, one iteration with
\eqref{eq:GSwave2} is the same as solving \eqref{eq:saddle_syst_wave}
with $K_{14}=0$.

\begin{figure}[htbp]
  % \psfrag{x}{$x$}
  % \psfrag{y}{$y$}
  % \centering
  % \includegraphics[height=0.47\textwidth]{pictures/wave2d/test1_mid_diff_mesh/avgCoeff.eps}
  % \includegraphics[height=0.47\textwidth]{pictures/wave2d/test1_mid_diff_mesh_noise/avgCoeff.eps}\\
  % \includegraphics[height=0.47\textwidth]{pictures/wave2d/test2_mid_diff_mesh/avgCoeff.eps}
  % \includegraphics[height=0.47\textwidth]{pictures/wave2d/test2_mid_diff_mesh_noise/avgCoeff.eps}\\
  \centering
  \includegraphics[width=1\textwidth]{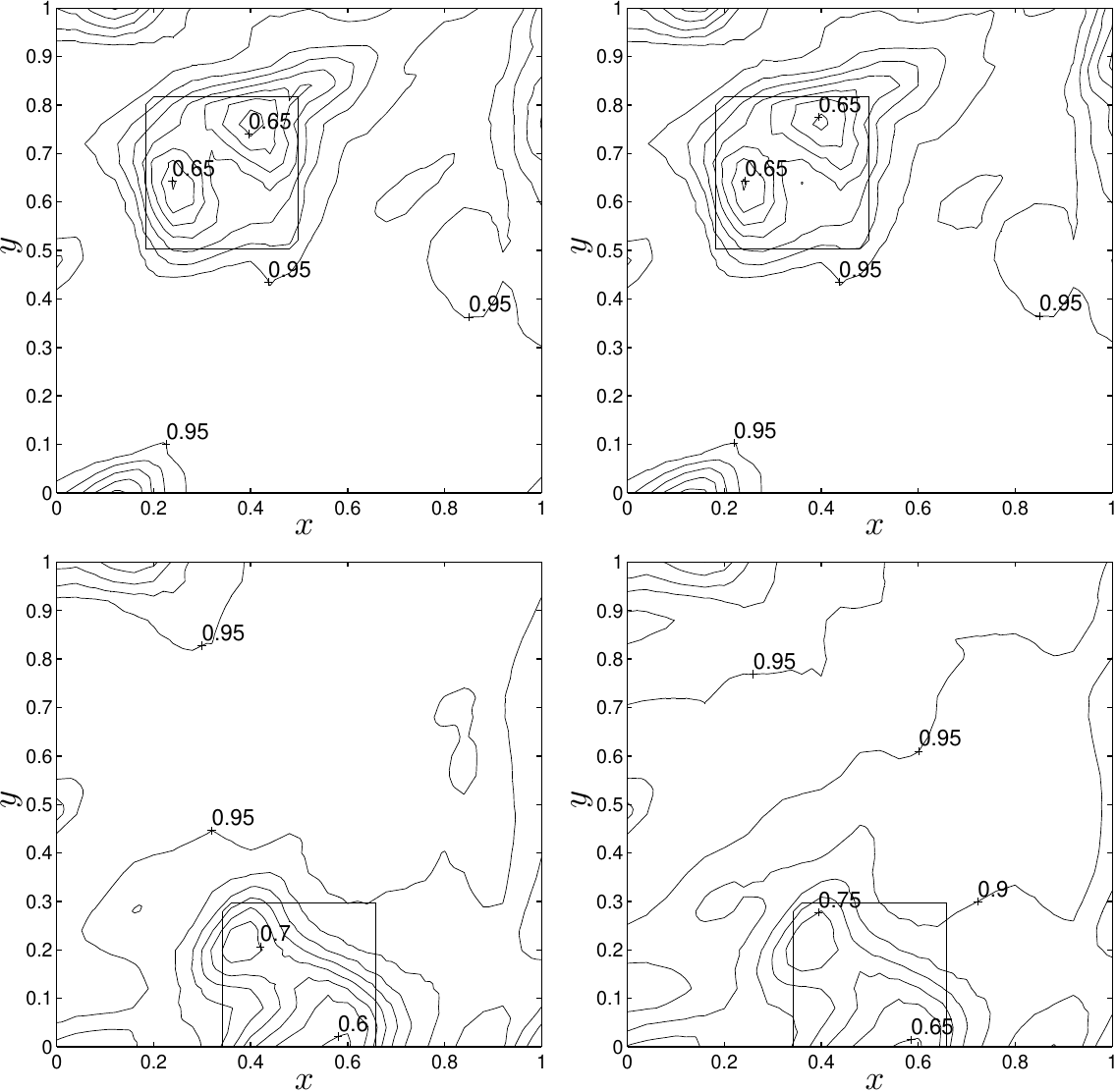}
  \caption{2D reconstruction using the weighted average
    \eqref{eq:avg3}, final time $T=1.5$ and
    Neumann boundary condition $2\sin( 4\pi t )$ for $x\in[0.4,0.6]$, $t<0.5$
    and $0$ elsewhere.
    The data $u^*$ was simulated by solving the forward equation on a
    quasi-uniform mesh with 3232 triangles and 328 time steps while the
    inverse problem was solved on a uniform mesh with 1250
    triangles and 30 time steps. 
    Measurements from the whole boundary were used.
    Top: Reconstruction of $\sigma_{true}=0.5$ inside the square
    $[0.2,0.5]\times[0.5,0.8]$ and $\sigma_{true}=1$ elsewhere, with no
    noise in data (left) and 10\% noise in data (right).
    Bottom: Reconstruction of $\sigma_{true}=0.5$ inside the square
    $[0.35,0.65]\times[0,0.3]$ and $\sigma_{true}=1$ elsewhere, with no
    noise in data (left) and 10\% noise in data (right).
  }
  \label{fig:2Dwave}
\end{figure}

\begin{figure}[htbp]
%   \psfrag{x}{$x$}
%   \psfrag{y}{$y$}
%   \centering
%   \includegraphics[height=0.25\textwidth]{pictures/wave2d/test1_mid_diff_mesh/measNew5.eps}
%   \includegraphics[height=0.25\textwidth]{pictures/wave2d/test1_mid_diff_mesh/measNew15.eps}
%   \includegraphics[height=0.25\textwidth]{pictures/wave2d/test1_mid_diff_mesh/measNew25.eps}\\
%   \includegraphics[height=0.25\textwidth]{pictures/wave2d/test1_mid_diff_mesh/coeff5.eps}
%   \includegraphics[height=0.25\textwidth]{pictures/wave2d/test1_mid_diff_mesh/coeff15.eps}
%   \includegraphics[height=0.25\textwidth]{pictures/wave2d/test1_mid_diff_mesh/coeff25.eps}
  \centering
  \includegraphics[width=1\textwidth]{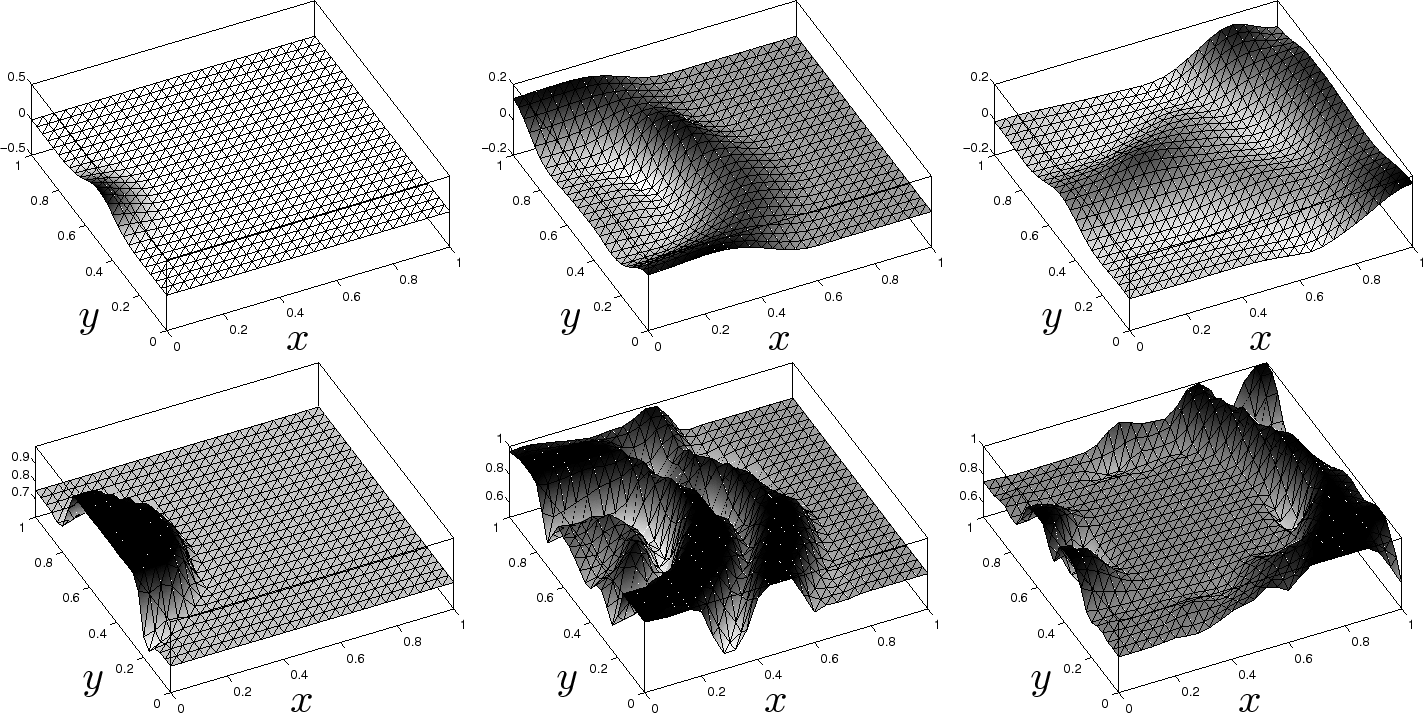}
  \caption{ Measurements $u^*$ (top) and $\mathfrak{h}_\delta'(\nabla u
    \cdot \nabla q)$ (bottom) for timesteps 5, 15 and 25.
    The data here corresponds to the top left plot in Figure
    \ref{fig:2Dwave}, and $u^*$ is interpolated onto the mesh used for the
    calculation of $u$ and $q$.
  }
  \label{fig:2Dwave2}
\end{figure}

In Figure \ref{fig:2Dwave}, a two dimensional example of reconstruction two
different speed coefficients is shown. The measured data was here simulated
by solving the wave equation for $\sigma_{true}$ with the symplectic
backward Euler method for \eqref{eq:constr_wave_syst}, which can be
written as the second order scheme 
\begin{equation*}
  \begin{aligned}
    \int_\Omega (u_{n+1}-2u_n+u_{n-1})w \ud x = 
    \int_{\partial\Omega} jw \ud s - \int_\Omega \sigma \nabla u_n
    \cdot \nabla w \ud x, \quad \forall w \in \bar V.
  \end{aligned}
\end{equation*}
Since the wave equation is a conservation law and is reversible in
time it is tempting to believe that it would be easier to control than
the heat equation but there are some computational drawbacks: numerical errors are
propagated in time and there seems to be many local minima. From the
approximation $\mathfrak{h}'_\delta(\nabla u\cdot \nabla q)$ in Figure \ref{fig:2Dwave2} it is evident
that the time dependent reconstruction varies a lot over time and is
not a good approximation of the time independent wave coefficient
$\sigma_{true}$.

%\section{Conclusions}

% Bibliography
\clearpage
\bibliographystyle{plain}

\begin{thebibliography}{10}

\bibitem{kunisch}
H.~T. Banks and K.~Kunisch.
\newblock {\em Estimation techniques for distributed parameter systems},
  volume~1 of {\em Systems \& Control: Foundations \& Applications}.
\newblock Birkh\"auser Boston Inc., Boston, MA, 1989.

\bibitem{barron}
Emmanuel~Nicholas Barron and Robert Jensen.
\newblock The {P}ontryagin maximum principle from dynamic programming and
  viscosity solutions to first-order partial differential equations.
\newblock {\em Trans. Amer. Math. Soc.}, 298(2):635--641, 1986.

\bibitem{bensoe}
M.~P. Bends{\o}e and O.~Sigmund.
\newblock {\em Topology optimization}.
\newblock Springer-Verlag, Berlin, 2003.
\newblock Theory, methods and applications.

\bibitem{saddle-point}
Michele Benzi, Gene~H. Golub, and J{\"o}rg Liesen.
\newblock Numerical solution of saddle point problems.
\newblock {\em Acta Numer.}, 14:1--137, 2005.

\bibitem{ghattas}
George Biros and Omar Ghattas.
\newblock Parallel {L}agrange-{N}ewton-{K}rylov-{S}chur methods for
  {PDE}-constrained optimization. {I}. {T}he {K}rylov-{S}chur solver.
\newblock {\em SIAM J. Sci. Comput.}, 27(2):687--713 (electronic), 2005.

\bibitem{borcea}
Liliana Borcea.
\newblock Electrical impedance tomography.
\newblock {\em Inverse Problems}, 18(6):R99--R136, 2002.

\bibitem{jesper}
Jesper Carlsson.
\newblock Pontryagin approximations for optimal design of elastic structures.
\newblock {\em preprint}, 2006.

\bibitem{css}
Jesper Carlsson, Mattias Sandberg, and Anders Szepessy.
\newblock Symplectic pontryagin approximations for optimal design.
\newblock Preprint.

\bibitem{cel}
M.~G. Crandall, L.~C. Evans, and P.-L. Lions.
\newblock Some properties of viscosity solutions of {H}amilton-{J}acobi
  equations.
\newblock {\em Trans. Amer. Math. Soc.}, 282(2):487--502, 1984.

\bibitem{engl}
Heinz~W. Engl, Martin Hanke, and Andreas Neubauer.
\newblock {\em Regularization of inverse problems}, volume 375 of {\em
  Mathematics and its Applications}.
\newblock Kluwer Academic Publishers Group, Dordrecht, 1996.

\bibitem{fenics}
{FE}ni{CS}.
\newblock {FE}ni{CS} project.
\newblock URL: url{http//www.fenics.org/}.

\bibitem{wanner}
Ernst Hairer, Christian Lubich, and Gerhard Wanner.
\newblock {\em Geometric numerical integration}, volume~31 of {\em Springer
  Series in Computational Mathematics}.
\newblock Springer-Verlag, Berlin, second edition, 2006.
\newblock Structure-preserving algorithms for ordinary differential equations.

\bibitem{dolfin}
J.~Hoffman, J.~Jansson, A.~Logg, and G.~N. Wells.
\newblock {DOLFIN}.
\newblock URL: url{http//www.fenics.org/dolfin/}.

\bibitem{lions}
J.-L. Lions.
\newblock {\em Optimal control of systems governed by partial differential
  equations.}
\newblock Translated from the French by S. K. Mitter. Die Grundlehren der
  mathematischen Wissenschaften, Band 170. Springer-Verlag, New York, 1971.

\bibitem{ss}
M.~Sandberg and A.~Szepessy.
\newblock Convergence rates of symplectic {P}ontryagin approximations in
  optimal control theory.
\newblock {\em M2AN}, 40(1), 2006.

\bibitem{mattias}
Mattias Sandberg.
\newblock Convergence rates for numerical approximations of an optimally
  controlled {G}inzburg-{L}andau equation.
\newblock {\em preprint}, 2006.

\bibitem{vogel}
Curtis~R. Vogel.
\newblock {\em Computational methods for inverse problems}, volume~23 of {\em
  Frontiers in Applied Mathematics}.
\newblock Society for Industrial and Applied Mathematics (SIAM), Philadelphia,
  PA, 2002.
\newblock With a foreword by H. T. Banks.

\end{thebibliography}

\end{document}